\documentclass[pdflatex,sn-mathphys-num]{sn-jnl}
\usepackage{graphicx} 
\usepackage{subcaption}
\usepackage[export]{adjustbox}
\usepackage{multirow}%
\usepackage{amsmath,amssymb,amsfonts}%
\usepackage{amsthm}%
\usepackage{mathrsfs}%
\usepackage[title]{appendix}%
\usepackage{xcolor}%
\usepackage{textcomp}%
\usepackage{manyfoot}%
\usepackage{booktabs}%
\usepackage{algorithm}%
\usepackage{algorithmicx}%
\usepackage{algpseudocode}%
\usepackage{listings}%
\usepackage{todonotes, comment}
\usepackage{lineno}
\usepackage{natbib} 
\theoremstyle{thmstyleone}%
\newtheorem*{theorem-non}{Theorem}

\theoremstyle{thmstyletwo}%
\newtheorem*{example}{Example}%

\newcommand{\scri}{\mathscr{I}}

\newcommand{\Op}{\operatorname{Op}}
\newcommand{\R}{\mathbb{R}}

\begin{document}


\title{From Penrose to Melrose: Computing Scattering Amplitudes at Infinity for Unbounded Media}

\author{\fnm{An\i l} \sur{Zengino\u{g}lu}}\email{anil@umd.edu}
\affil{\orgdiv{Institute for Physical Science \& Technology}, \orgname{University of Maryland}, \orgaddress{\street{4254 Stadium Drive}, \city{College Park}, \postcode{20742}, \state{MD}, \country{USA}}}

\abstract{
We develop a method to compute scattering amplitudes for the Helmholtz equation in variable, unbounded media with possibly long-range asymptotics. Combining Penrose's conformal compactification and Melrose's geometric scattering theory, we formulate the time-harmonic scattering problem on a compactified manifold with boundary and construct a two-step solver for scattering amplitudes at infinity. The construction is asymptotic: it treats a neighborhood of infinity, and is meant to couple to interior solvers via domain decomposition. The method provides far-field data without relying on explicit solutions or Green's function representation. Scattering in variable media is treated in a unified framework where both the incident and scattered fields solve the same background Helmholtz operator. Numerical experiments for constant, short-range, and long-range media with single-mode and Gaussian beam incidence demonstrate spectral convergence of the computed scattering amplitudes in all cases.
}



\maketitle


\pagebreak


\section{Introduction}\label{sec:intro}

We study time-harmonic scattering by an obstacle embedded in a variable medium. Let $d\ge2$. A Helmholtz scattering problem consists of a bounded Lipschitz domain $\mathcal D\subset\mathbb R^{d}$ with boundary $\Gamma=\partial\mathcal D$ and exterior domain
$\mathcal D^{c}=\mathbb R^{d}\setminus\overline{\mathcal D}$, a real-valued
refractive index $n:\mathbb R^{d}\to\mathbb R$, and a boundary condition on $\Gamma$ (Dirichlet, Neumann, or impedance). A scattering solution at frequency $k>0$ is a function
$U:\mathcal D^{c}\to\mathbb C$ satisfying 
\begin{equation}\label{eq:helmholtz}
(\Delta + k^2 n(x)^2)\,U = 0 
\quad \text{in } \mathcal{D}^c,
\end{equation}
with the boundary condition on $\Gamma$, and an appropriate radiation condition at infinity. The asymptotic behavior of solutions depends on the refractive index. In this paper, we consider asymptotically homogeneous media with a refractive index satisfying
\begin{equation}\label{eq:refractive}
n(x)^2=1+b(x), \qquad b(x)\to 0 \quad\text{as } r=|x|\to\infty,
\end{equation}
uniformly in the angular variable $\omega=x/|x|\in\mathbb S^{d-1}$. 

Sommerfeld showed in 1912 that, unlike typical elliptic problems, the Helmholtz problem has non-unique solutions even when the field vanishes at infinity  \cite{sommerfeld1912greensche, schot1992eighty}. To ensure uniqueness, one must supplement the equation with a radiation condition that excludes waves arriving from infinity. Sommerfeld argued that physically acceptable solutions are purely outgoing at infinity. This principle has guided a century of analysis and computation in scattering theory and numerics.
However, in modern scattering theory, incoming waves from infinity are not considered unphysical. On the contrary, the mathematically natural data are \emph{prescribed} at infinity: one specifies the incoming input and reads off the outgoing output via the \emph{scattering map} on the celestial sphere at infinity \cite{melrose1995geometric, dyatlov2019mathematical}.

Despite the central role of this structure in modern scattering theory, its computational use is limited. The conventional frequency-domain approach models incidence by prescribing a global, free-space solution that induces boundary data on the surface of the obstacle. However, in unbounded media with long-range modifications, the incident solution is not a priori given. In this paper, we numerically construct the incident field from incoming radiation data as a first step in the computation of scattering amplitudes for long-range media.

\begin{figure}[t]
  \centering
  \begin{subfigure}[t]{0.44\textwidth}
    \centering
    \includegraphics[width=\linewidth, valign=c]{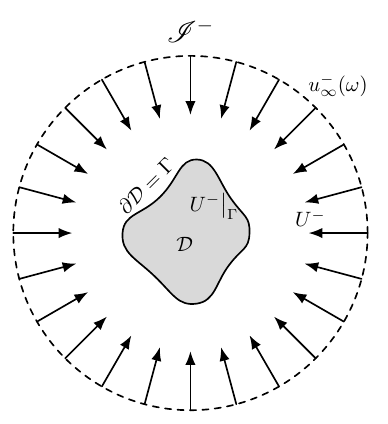}
    \caption{Incoming radiation from past null infinity, $u_\infty^-(\theta)$, induces an incident field at the obstacle boundary, $U^-\big|_\Gamma$.}
    \label{fig:incident}
  \end{subfigure}\hfill
  \begin{subfigure}[t]{0.44\textwidth}
    \centering
    \includegraphics[width=\linewidth, valign=c]{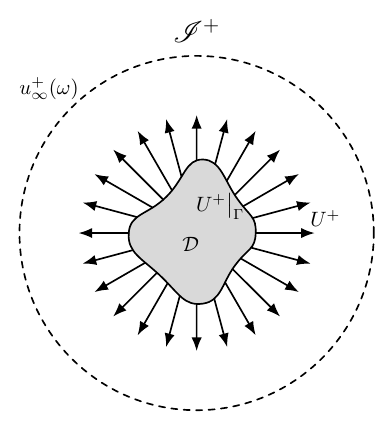}
    \caption{Outgoing radiation from the obstacle, $U^+\big|_\Gamma$, induces the scattered far-field at future null infinity $u_\infty^+(\theta)$.}
    \label{fig:scattered}
  \end{subfigure}
  \caption{We solve the time-harmonic scattering problem \eqref{eq:helmholtz}, where incoming radiation from infinity, $u^-_\infty(\omega)$ is scattered by an obstacle $\mathcal{D}$ and is measured at infinity as the far field $u^+_\infty(\omega)$. The computation consists of two steps based on the separation of the total field $U$ into incident $(-)$ and scattered $(+)$ fields: $U=U^-+U^+$. In the first step, we compute the incident field $U^-$ by prescribing incoming boundary conditions at past null infinity, $\mathscr{I}^{-}$, and extract the data on the obstacle boundary $\Gamma=\partial\mathcal{D}$. In the second step, we use the obstacle boundary data to solve for the scattered field $U^+$ and read off the far field at future null infinity, $\mathscr{I}^{+}$. The singularity of the map between $U^\pm$ and $u^\pm_\infty$ is absorbed by a suitable weight function \eqref{eq:scale_out} that corresponds to regular compactification at null infinity $\mathscr I^\pm$ described in Sec.~\ref{sec:penrose}.}
  \label{fig:incident-scattered}
\end{figure}

A key observation from geometric scattering theory is that \emph{radiation conditions are local at infinity}, in the sense that a clean split of incoming and outgoing radiation is only valid at infinity \cite{melrose1995geometric, dyatlov2019mathematical, epstein2024solving}. To make this idea concrete, assume that \(U\) solves the constant-index Helmholtz equation in the exterior of the scatterer, i.e.~Eq.~\eqref{eq:helmholtz} with $n(x)=1$. At large radius \(r=|x|\), the general solution admits an incoming/outgoing decomposition,
\begin{equation}\label{eq:decomposition}
U(r,\omega) \sim r^{-\frac{d-1}{2}}\left( e^{-ikr}\,u_\infty^-(\omega)+e^{+ikr}\,u_\infty^+(\omega)\right),
\qquad \text{as } \ r\to\infty.
\end{equation}
The scattering map \(S(k)\) is the operator mapping the incoming incident data \(u_\infty^-(\omega)\) to the outgoing scattered solution \(u_\infty^+(\omega)=S(k)u_\infty^-(\omega)\). A local decomposition of a regular field into pointwise incoming/outgoing parts is not available in general \cite{felbacq2025local, martin2018out}. In the microlocal framework, the mathematically natural datum is prescribed at infinity as the incoming profile $u_\infty^-(\omega)$. The computational framework implied by the mathematical analysis of the scattering problem is depicted as a two-step solver in Fig.~\ref{fig:incident-scattered}. We assume a split of the total field into an incident and scattered part as $U=U^-+U^+$. Numerically, we prescribe $u_\infty^-(\omega)$ at $\mathbb{S}^{d-1}_\infty$, solve for the corresponding incident field $U^-$ in the background medium, impose the physical boundary condition on the scatterer to solve for $U^+$, and read off $u_\infty^+(\omega)$ at $\mathbb{S}^{d-1}_\infty$. This picture aligns the computational approach with the continuous scattering map $u_\infty^-\mapsto u_\infty^+$.

An immediate difficulty with this procedure is how to handle incoming data at infinity. Particularly for long-range media, the incident field can not be prescribed as an undisturbed wave of the free-space equation; one needs to solve the variable-coefficient operator on an unbounded domain. In geometric scattering \cite{melrose1995geometric}, one introduces scattering metrics by radial compactification of the form $r=1/\rho $. This transformation leads to a singular operator at the boundary. For demonstration, consider the Helmholtz equation written in spherical coordinates $(r,\omega)$ as
\[
\partial_r^2 U + \frac{d-1}{r}\partial_r U + \frac{1}{r^2}\,\Delta_{\mathbb S^{d-1}} U + k^2 U = 0.
\]
Applying radial compactification and dividing the equation by $\rho^2$, we obtain a singular operator,
\begin{equation}\label{eq:singular-helmholtz}
(\rho\,\partial_\rho)^2 U +(2-d)(\rho\,\partial_\rho) U +\Delta_{\mathbb S^{d-1}} U
+\frac{k^2}{\rho^2} U =0.
\end{equation}
The effective wave number, $k/\rho$, is singular because time-harmonic solutions oscillate infinitely many times on an infinite domain. Radial compactification translates the unbounded domain problem to an unbounded resolution problem, and is therefore not suitable for direct numerical implementation. Mathematically, the singularity structure of the above equation is treated by scattering calculus (see App.~\ref{sec:spacetime}). We show that the ideas underlying scattering calculus can be implemented numerically by exploiting the connection of the formalism to the global causal structure of spacetimes.

Infinity has a richer structure in spacetime than in space, as revealed by Penrose \cite{penrose_asymptotic_1963, penrose2011republication}. We distinguish between infinity in space, time, and the characteristic directions. Infinity along the characteristic (or null) directions is called \textbf{null infinity}, denoted by \(\mathscr I\) (read ``scri'' for ``script i"). Null infinity has two components: past null infinity, $\mathscr I^-$, where incoming characteristics come from, and future null infinity, $\mathscr I^+$, where outgoing characteristics go to. Past null infinity is the natural location to prescribe incoming radiation data, while future null infinity is the natural location to read off outgoing radiation data.

Recently, techniques from relativity exploiting spacelike approach to null infinity (see \cite{zenginouglu2008hyperboloidal, zenginouglu2011geometric, zenginouglu2011hyperboloidal, jaramillo2021pseudospectrum, assaad2025quasinormal, panosso2024hyperboloidal, panosso2025hyperboloidal}) have been adapted to Helmholtz equations, allowing the computation of scattered solutions at future null infinity through \emph{null infinity compactification} (NIC) \cite{zenginouglu2021null}. The procedure involves scaling out the oscillatory decay along the radial direction\footnote{In general, $h$ can depend on the angles $\omega$. We consider the simpler radial case in this paper.} via
\begin{equation}\label{eq:scale_out} U^\pm(r,\omega) = \frac{e^{i k\,h(r)}}{r^{\frac{d-1}{2}}} u^\pm(r,\omega), \qquad h'(r) \to \pm 1,\end{equation}
where $h$ is called the height function (see Sec.~\ref{sec:penrose}). After the rescaling, the radially compactified Helmholtz equations for $u^\pm$ extend smoothly to the boundary. For example, with the choice of $h(r)=\pm r$, we obtain
\[ \rho^2\partial_\rho^2 u^\pm
+ 2\left(\rho\mp ik\right)\partial_\rho u^\pm
+ \Delta_{\mathbb S^{d-1}} u^\pm
+ \frac{(d-1)(3-d)}{4} u^\pm = 0
\]
The radiation condition becomes a local boundary condition, allowing us to \emph{make infinity part of the computational formalism}. In particular, NIC allows direct evaluation of the asymptotic quantities under suitable choices of the free functions such that $\lim_{r\to\pm\infty} u^\pm(r,\omega) = u^\pm_\infty(\omega)$. In this paper, we use this method to compute scattering amplitudes for obstacles in variable, unbounded media with incident data prescribed at infinity. The numerical experiments focus on two spatial dimensions but the method applies to arbitrary dimensions.

\textbf{Main contributions.} 

\begin{enumerate}
  \item We exploit a connection between Penrose's conformal compactification and Melrose's scattering calculus to construct a two-step solver for Helmholtz scattering on unbounded domains (see Fig.~\ref{fig:incident-scattered}). The solver relies on scaling out the oscillatory decay from asymptotic solutions using \eqref{eq:scale_out} combined with radial compactification.

  \item We compute scattering amplitudes from incoming radiation prescribed at infinity as a smooth profile with finite energy, \(u_\infty^-(\theta)\in L^2(\mathbb S^1)\). The far field \(u_\infty^+(\theta)\) is delivered as the output, providing the scattering map \(u_\infty^- \mapsto u_\infty^+=S(k)u_\infty^-\).
  
  \item We solve scattering problems for obstacles in variable, unbounded media, including short-range and Coulomb-type long-range cases. The numerical solver recovers the scattering amplitudes at infinity with spectral accuracy in all cases for single-mode and Gaussian beam incidence.
\end{enumerate}

The numerical calculations on unbounded domains rely on the following
\begin{theorem-non} Let $d\ge2$ and let $U:\mathcal D^{c}\to\mathbb C$ be a sufficiently regular solution of the variable-coefficient Helmholtz equation \eqref{eq:helmholtz} with real-valued refractive index $n=n(x)$ given at \eqref{eq:refractive} and augmented with appropriate boundary conditions on $\Gamma$ and a radiation condition at infinity.

We let $r=r(\rho):(0,1)\to(0,\infty)$ be a $C^2$ strictly increasing radial coordinate such that
\begin{equation}\label{eq:compactification}
  r(0)=0, \qquad
  \lim_{\rho\to1^-} r(\rho)=\infty,
  \qquad
\lim_{\rho\to1^-}\frac{r^2(\rho)}{r'(\rho)}=\eta
\end{equation}
for some fixed $\eta>0$. 
Further, let $h=h(r):(0,\infty)\mapsto \mathbb R$ be a $C^2$ radial height function. Denote $H(\rho):= h'(r(\rho))$.
Express $n(x)$ in terms of the radial variables $n(\rho,\omega):=n(r(\rho)\omega)$, $(r,\omega) \in \mathbb{R}^+\times \mathbb{S}^{d-1}$ and denote $G(\rho):=\frac{1}{r'(\rho)}$. Assume the following regularity condition holds
\begin{equation}\label{eq:rate-condition}
n^2(\rho,\omega)-H^2(\rho) = \mathcal O\big(G(\rho)\big)
\quad\text{as }\rho\to1,
\end{equation}
uniformly in $\omega\in\mathbb S^{d-1}$ (note that $G(\rho)\stackrel{\rho\to 1^-}{\longrightarrow} 0$). Finally, define the rescaled field
\begin{equation}\label{eq:weight-rescaling-NIC}
u(\rho,\omega) = \frac{U(r(\rho)\omega)}{w(r(\rho))} \quad \text{with} \quad w(r(\rho)) := \exp\left(ik h(r(\rho))\right)\,r^{-\frac{d-1}{2}}(\rho).
\end{equation}
Then the Helmholtz equation \eqref{eq:helmholtz} with \eqref{eq:refractive} is expressed in terms of the rescaled field $u(\rho,\omega)$ as 
\begin{align}
\label{eq:NIC-Helmholtz}
& \partial_{\rho}\big(G(\rho)\,\partial_{\rho}u\big)
   + 2ik\,H(\rho)\,\partial_{\rho}u
   + \frac{1}{G(\rho)\,r^2(\rho)}\Delta_{\mathbb S^{d-1}}u \nonumber \\
&\qquad
   + \left[
      k^2\,\frac{n^2(\rho,\omega)-H^2(\rho)}{G(\rho)}
      + ik\,H'(\rho)
      + \frac{(d-1)(3-d)}{4\,G(\rho)\,r^2(\rho)}
     \right]u = 0. 
\end{align}
\end{theorem-non}

\noindent \textbf{Remarks:} The reasoning behind the conditions in the theorem is explained in Sec.~\ref{sec:spacetime}. Examples for radial compactifications satisfying \eqref{eq:compactification} are \eqref{eq:compactify} and \eqref{eq:hyperbolic-radius}. Examples for height functions satisfying \eqref{eq:rate-condition} are $h(r)=r$ and \eqref{eq:hyperbolic-height}. The assumption \eqref{eq:compactification} dictates the rate at which $G(\rho)$ approaches zero, ensuring that $1/r(\rho)$ is an asymptotic boundary defining function. The assumption \eqref{eq:rate-condition} dictates the rate at which $n^2(\rho,\omega)-H^2(\rho)$ approaches zero, ensuring that time surfaces approach the null cone.
Under these assumptions, all coefficients extend continuously to $\rho=1$. The formally second-order term in the $\rho$ direction reduces to first order because $G(1)=0$.

\section{Null infinity compactification}\label{sec:nic}
In this section, we discuss the computational framework of null infinity compactification (NIC) that leads to the NIC-Helmholtz equation \eqref{eq:NIC-Helmholtz}. We solve the scattering problem by first solving for the incident field with incoming data prescribed at infinity and then solving for the scattered field with boundary data at the obstacle. We contrast NIC with the standard approach for the constant and variable media cases because they solve different but related problems.

\subsection{The transformed Helmholtz equation}
Consider the Helmholtz equation with a (possibly variable) refractive index,
\begin{equation}\label{eq:helmholtz-d}
\Delta U(x) + k^2 n^2(x)\,U(x) = 0, \qquad x\in \mathbb{R}^d.
\end{equation}
We assume that the refractive index is real-valued and tends to a constant in every direction \eqref{eq:refractive}. Writing $x = r\omega$ with $r = |x| > 0$ and $\omega \in S^{d-1}$, we regard
$n$ as a function $n(r\omega)$ and assume that it approaches a constant as $r\to\infty$ uniformly for all $\omega \in S^{d-1}$. Taking this constant to be unity without loss of generality, we write
\begin{equation}\label{eqn:refractive_index}
n^2(x) = 1 + b(x), \qquad b(x) \to 0 \ \text{as } r\to\infty.
\end{equation}
In spherical coordinates $(r,\omega)$, the Helmholtz equation can be written as,
\begin{equation}
\partial_r^2 U + \frac{d-1}{r}\partial_r U + \frac{1}{r^2}\,\Delta_{\mathbb S^{d-1}} U 
+ k^2 n^2 U = 0.
\end{equation}
By \eqref{eqn:refractive_index}, there is no preferred direction at infinity and the dominant oscillation is governed by a radial phase. To capture the asymptotic behavior of the solution by the operator, we scale out the oscillatory decay with a radial weight \cite{zenginouglu2011geometric, zenginouglu2021null},
\begin{equation}\label{eq:w-rescale-d}
U(r,\omega) = w(r)\,u(r,\omega), 
\qquad 
w(r) := \frac{e^{\,ik\,h(r)}}{r^{\frac{d-1}{2}}}, 
\qquad 
H(r) := h'(r).
\end{equation}
In constant media, the choice $h(r)=r$ corresponds to the basis functions of the infinite element method. However, here we are transforming the Helmholtz operator. By direct computation, we get (compare \eqref{eq:Lhk_rad})
\[
\partial_r^2 u+ 2 i k H \,\partial_r u+ \frac{1}{r^2} \Delta_{\mathbb{S}^{d-1}} u + 
\left( k^2 \left(n^2-H^2\right) + i k H_{r} + \frac{(d-1)(3-d)}{4 r^2} \right) u= 0.
\]
Compactifying the spatial domain using \eqref{eq:compactification} gives \eqref{eq:NIC-Helmholtz}. For the numerical calculations in this paper, we restrict to two dimensions, \(d=2\), and radial refractive indices $n(x)=n(r)$. Except for the example in Sec.~\ref{sec:hyperbolic}, we choose the radial compactification as 
\begin{equation}\label{eq:compactify}
r=\frac{\rho}{1-\rho},\qquad 
G(\rho):=\frac{d\rho}{dr}= (1-\rho)^2,
\end{equation}
To satisfy the regularity condition \eqref{eq:rate-condition}, we choose the boost function as \(H = \pm n\). Then, the \(k^2\) term vanishes and we get
\begin{equation}\label{eq:u-eq-2D-compact}
\partial_\rho\left( (1-\rho)^2 \partial_\rho u \right) \pm 2ik n\,\partial_\rho u + \frac{1}{\rho^2} \partial_\theta^2 u + \left(\pm ik n_\rho + \frac{1}{4\rho^2}\right) u = 0.
\end{equation}
The positive sign corresponds to outgoing waves to future null infinity; the negative sign to incoming waves from past null infinity. For the incoming choice, boundary data at \(\rho=1\) model incident waves from infinity, while for the outgoing case, data are imposed only at the obstacle boundary. In the following, we treat the constant index and variable index cases separately, as there are significant differences in the scattering behavior and the standard treatment.

\subsection{Constant media}\label{sec:2D-constant}
We first consider the constant-medium problem with \(n\equiv 1\). 
We assume an obstacle with a sound-soft boundary. The total field \(U\) satisfies the homogeneous Helmholtz equation in the exterior domain \(\mathcal{D}^c\) and the Dirichlet boundary condition on the obstacle surface \(\Gamma\):
\[
(\Delta + k^2) U = 0 \quad \text{in } \mathcal{D}^c,
\qquad
U = 0 \quad \text{on } \Gamma.
\]
The solution satisfies the following asymptotic expansion near infinity,
\begin{equation}
\label{eq:radiation-expansion}
U(r,\theta) = r^{-1/2}\left(e^{-ik r} u^-_\infty(\theta)+ e^{+ik r} u^+_\infty(\theta)\right) +\mathcal{O}\left(r^{-3/2}\right), \qquad r\to\infty.
\end{equation}
The amplitudes $u^\pm_\infty\in L^2(\mathbb{S}^1)$ are the outgoing/incoming radiation data. By linearity of the Helmholtz equation, we decompose the total field $U$ into an incident field and a scattered field, \(U = U^- + U^+\). At this point, we present the standard setup and the two-step solver separately.

\subsubsection{Standard setup for constant media}\label{sec:standard}
In the standard setup, \(U^-\) is prescribed a priori as a free-space Helmholtz solution,
\begin{equation}\label{eq:free-space}
(\Delta + k^2) U^- = 0 \quad \text{in } \mathbb{R}^{d}
\qquad \text{(e.g. a plane wave \(U^-(r,\theta)=e^{ik r\,\cos(\theta-\theta_0)}\)).}
\end{equation}
The scattered field \(U^+\) then solves a homogeneous exterior boundary-value problem with Dirichlet data induced by \(U^-\) on the obstacle surface and satisfies the Sommerfeld radiation condition at infinity:
\[
\begin{cases}
(\Delta + k^2)U^+ = 0 & \text{in } \mathcal{D}^c,\\[2pt]
U^+ = -\,U^- & \text{on } \Gamma,\\[4pt]
\sqrt{r}\,\left(\partial_r - i k\right)\,U^+ = 0 & \text{as } r\to\infty.
\end{cases}
\]
Under the radiation condition and standard smoothness assumptions, \(U^+\) admits the far-field expansion
\[
U^+(x) = \frac{e^{ik r}}{\sqrt{r}}\,u^+_\infty(\theta)
+ \mathcal{O}\left(r^{-3/2}\right),
\qquad r\to\infty,
\]
defining the far-field pattern \(u_\infty^+:\mathbb{S}^1\to\mathbb{C}\).
The total field is recovered as \(U = U^- + U^+\). Note that both $U^+$ and $U^-$ solve the same free-space Helmholtz operator.

\subsubsection{NIC setup for constant media}
To compute the incident and scattered fields directly on the infinite domain, we define the rescaled variables \(u^\pm\) via \eqref{eq:w-rescale-d} with $h=\pm r$ and $d=2$, 
\begin{equation}\label{eq:u-plus-minus}
U^- = w^- u^- = \frac{e^{-ik r}}{\sqrt{r}}\,u^-,\qquad U^+ = w^+ u^+ = \frac{e^{+ik r}}{\sqrt{r}}\,u^+.
\end{equation}
We consider the rescaled variables as functions of the compactifying coordinate $\rho$. Then, radiation data $u^\pm_\infty(\theta)$ in the asymptotic expansion \eqref{eq:radiation-expansion} are obtained by direct evaluation at infinity: $u^\pm(1,\theta)=u^\pm_\infty(\theta)$. Using $n=1$ in \eqref{eq:u-eq-2D-compact}, we get the equation,
\begin{equation}\label{eq:compactified}
\partial_\rho \left( (1-\rho)^2 \partial_\rho u^\pm\right) \pm 2ik\,\partial_\rho u^\pm + \frac{1}{\rho^2} \partial_\theta^2 u^\pm + \frac{1}{4\rho^2}u^\pm=0.
\end{equation}
The operator becomes first order in $\rho$ at \(\rho=1\) with the boundary compatibility condition
\begin{equation}\label{eq:compat-2D}
\pm\,2ik\,\partial_\rho u^\pm(1,\theta) + \partial_\theta^2 u^\pm(1,\theta) + \frac{1}{4}u^\pm(1,\theta)=0.
\end{equation}
For the incoming choice with the negative sign in the equation, we prescribe the amplitude as
\begin{equation}\label{eq:incoming-trace}
u^-(1,\theta)=u_\infty^-(\theta)\in L^2(\mathbb S^1),
\end{equation}
and impose no further condition at \(\rho=1\); the equation supplies \(u^-_\rho(1,\theta)\) via \eqref{eq:compat-2D}.
We solve \eqref{eq:compactified} for $u^-$ on a punctured domain. We read off the incident field on the surface representing the obstacle boundary, $\Gamma$, as \(U^-|_\Gamma = w^- u^- |_\Gamma\) where $w^\pm = e^{\pm i k r}/\sqrt{r}$. This incident field is then used to solve the scattering problem for the scattered field \(U^+\). In the sound-soft case, we impose the physical boundary condition for the total field $U$ on the obstacle surface \(\Gamma\) as \(U|_\Gamma=0\), which implies a relationship between the incident and scattered fields as $U^+|_\Gamma = -U^-|_\Gamma $. The boundary data for the rescaled scattered NIC field is obtained by
\begin{equation}\label{eq:numerical_boundary_data} 
  u^+|_{\Gamma} = (w^+)^{-1} U^+|_\Gamma = - (w^+)^{-1} U^-|_\Gamma =  - (w^+)^{-1} w^- u^-|_\Gamma = -e^{-2 i k r} u^-|_\Gamma. 
\end{equation}
The outgoing radiation data at infinity, and hence the scattering amplitudes, are then read off simply by evaluating $u^+$ at null infinity,
\begin{equation}\label{eq:far-field}
u^+(1,\theta)=u_\infty^+(\theta)\in L^2(\mathbb S^1),
\end{equation}
This step completes the solution of the scattering problem $u_\infty^-\mapsto u_\infty^+$ as depicted in Fig.~\ref{fig:incident-scattered} mapping  \eqref{eq:incoming-trace} to  \eqref{eq:far-field}. This two-step process reflects the foundational structure of the continuous theory described in Sec.~\ref{sec:melrose}.
However, in constant media, computing the incident field numerically is unnecessary because it can always be constructed analytically. This computation is performed in Sec.~\ref{sec:single-mode-constant} primarily as a testbed for the two-step solver.

\subsection{Variable media}\label{sec:2D-variable}
The geometric setup for variable media is the same as in Sec.~\ref{sec:2D-constant} and Fig.~\ref{fig:incident-scattered} but with a non-vanishing perturbation $b$ in \eqref{eqn:refractive_index}. As in the constant-medium case, the total field \(U\) satisfies the variable-coefficient Helmholtz equation in the exterior domain $\mathcal{D}^c$ with Dirichlet boundary condition on \(\Gamma\):
\begin{equation}\label{eq:variable-helmholtz}
(\Delta + k^2 n^2(x)) U = 0 \quad \text{in } \mathcal{D}^c,
\qquad
U = 0 \quad \text{on } \Gamma.
\end{equation}
The asymptotic expansion \eqref{eq:radiation-expansion} is valid for short-range refractive indices. For long-range cases where the perturbation decays slowly as $r$, the phases in \eqref{eq:radiation-expansion} acquire a Coulomb-type, logarithmic modification. Next, we decompose the total field as before into incident and scattered parts, \(U=U^-+U^+\), and describe the standard and NIC setups separately. 

\subsubsection{Standard setup for variable media}\label{sec:standard_variable}
In the standard setup, \(U^-\) is assumed to satisfy the \emph{free-space} Helmholtz equation \eqref{eq:free-space} even though the background medium has a variable refractive index. This assumption implies a source term for the scattered field \(U^+\) in the exterior domain:
\begin{equation} \label{eq:sourced-helmholtz}
(\Delta + k^2(1+b))\,U^+ = -k^2\,b\,U^- \quad \text{in } \mathcal{D}^c.
\end{equation}
This equation is known as the Lippmann--Schwinger equation \cite{colton2013inverse}. The boundary condition on the obstacle surface is derived from the sound-soft condition as before: $U^+|_\Gamma=-U^-|_\Gamma$. The radiation condition depends on the asymptotics of \(b\). For a short-range background with \(b(x)=\mathcal O(r^{-1-\varepsilon})\), we impose the usual Sommerfeld condition:
\[
\sqrt r\,(\partial_r - i k)\,U^+ \longrightarrow 0 \quad \text{as } r\to\infty.
\]
For a long-range background with a \(1/r\) tail such as
\(
n^2(r)=1+\frac{\kappa}{r},
\)
implying $n(r)=1+\frac{\kappa}{2r}+\mathcal{O}(r^{-2})$, the outgoing condition is Coulomb-modified:
\begin{equation}\label{eq:longSommerfeld}
\sqrt r\left[\partial_r - i k \left(1+ \frac{\kappa}{2 r}\right)\right]\,U^+ \longrightarrow 0
\quad \text{as } r\to\infty.
\end{equation}
While the total field produced by the above formulation satisfies \eqref{eq:variable-helmholtz}, the decomposition does not describe propagation in the actual medium: \(U^-\) solving \eqref{eq:free-space} is not a solution of the variable-index Helmholtz equation, hence the boundary datum \(U^-|_\Gamma\) is not the trace of a field that has propagated through the true background. The Lippmann--Schwinger formulation puts all scattering, including backscatter from the medium, into the scattered field \(U^+\). In contrast, geometric scattering theory treats incident and scattered fields on equal footing as solutions of the same Helmholtz operator. Consequently, we require that both components solve the same equation in the exterior,
\[
(\Delta + k^2 n^2(x))\,U^- = 0,\qquad
(\Delta + k^2 n^2(x))\,U^+ = 0 \quad \text{in } \mathcal{D}^c.
\]
Solving for both fields numerically is generally unnecessary for scattering problems in short-range media. However, for long-range media where the physical incoming wave carries a logarithmic phase correction, an ad hoc solution of the free-space Helmholtz equation should not be used globally. Imposing the corrected outgoing condition \eqref{eq:longSommerfeld} for \(U^+\) fixes only the outgoing part; a free-space incident solution does not have the correct phase. A consistent formulation for long-range media requires the computation of the incident and scattered fields using the same operator or the construction of a global incident field that satisfies the modified long-range asymptotics.


\subsubsection{NIC setup for variable media}\label{sec:nic_variable}
Let $\mathcal{N}$ be a function with $\mathcal{N}'(r)=n(r)$. Choose the height function as $h(r)=\pm \mathcal{N}(r)$. This choice ensures that the compactification is performed along the characteristics. The regularity condition \eqref{eq:rate-condition} is trivially satisfied. The rescaled fields are defined as
\begin{equation}\label{eq:u-variable-rescaled}
U^- = \frac{e^{-ik \mathcal{N}(r) }}{\sqrt{r}}\,u^-,\qquad U^+ = \frac{e^{+ik \mathcal{N}(r)}}{\sqrt{r}}\,u^+, \quad \text{where} \quad \mathcal{N}(r):=\int n(r)dr.
\end{equation}
The compactified equation reads
\begin{equation}\label{eq:compactified-variable}
\partial_\rho \left( (1-\rho)^2 \partial_\rho u^\pm\right) \pm  2ik n\,\partial_\rho u^\pm + \frac{1}{\rho^2}\partial_\theta^2 u^\pm + \left(\pm i k n_\rho + \frac{1}{4\rho^2}\right)u^\pm=0.
\end{equation}
The radiation condition is naturally part of the compactified equation and arises as a boundary compatibility condition,
\begin{equation}\label{eq:compat-2D-variable}
\pm\,2ik n(1,\theta)\,\partial_\rho u^\pm(1,\theta) + \partial_\theta^2 u^\pm(1,\theta) + \left(\pm i k n_\rho(1,\theta) + \frac{1}{4} \right) u^\pm(1,\theta)=0.
\end{equation}
The boundary data for the scattered field is obtained from the incident field by 
\begin{equation}\label{eq:boundary-variable}
  u^+|_\Gamma=-e^{-2ik\mathcal{N}(r)} u^-|_\Gamma.
\end{equation}

\pagebreak

\section{Numerical experiments}\label{sec:numerics}
\begin{figure}[t]
  \centering
  \includegraphics[width=\textwidth]{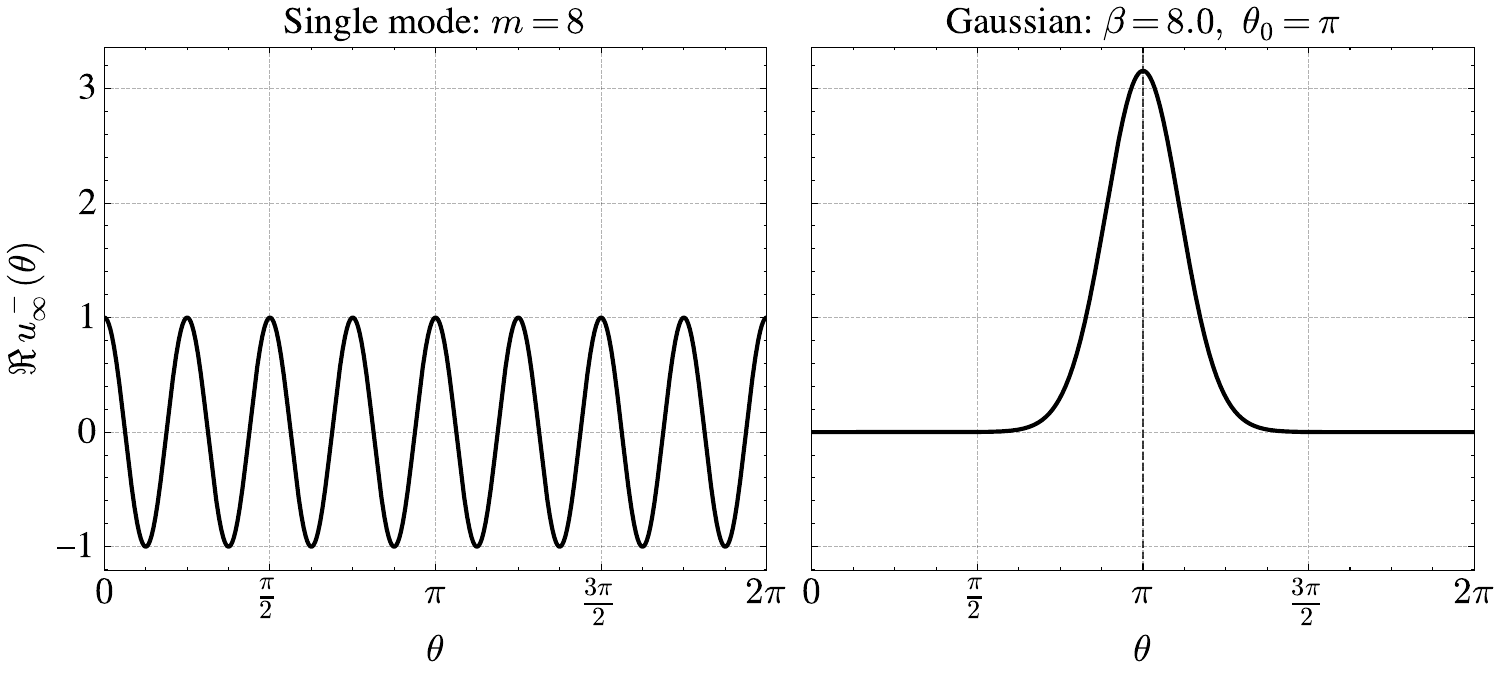}
  \caption{Real parts on a shared axis for the two types of incoming data, $\Re u_\infty^-(\theta)$, prescribed at past null infinity \(\mathscr I^-\): single mode \eqref{eq:single-mode} with $m=8$ and a circular Gaussian beam using a von Mises distribution \eqref{eq:gaussian} with $\beta=8, \ \theta_0=\pi$.}
  \label{fig:incoming-panel}
\end{figure}

We demonstrate the numerical solution of the NIC-transformed Helmholtz equation with incoming radiation from infinity on two families of incident data: single Fourier modes, and circular Gaussian beams (Fig.~\ref{fig:incoming-panel}). 
We prescribe incident radiation at the celestial circle\footnote{Our experiments are in two spatial dimensions, but application to three spatial dimensions with data prescribed at the celestial sphere is straightforward.} at null infinity. 
To compare the numerical solutions to analytic ones, we restrict our calculations to scattering at the unit disk and radial media. The single-mode computations help us isolate the centrifugal barrier and demonstrate the effect of variable media on the scattering map. For more realistic data at infinity, we also compute scattering amplitudes with circular Gaussian beam incidence in Sec.~\ref{sec:gaussian}.

\subsection{Single mode}\label{sec:single-mode}
The simplest smooth choice of incoming data from infinity is a single Fourier mode,
\begin{equation}\label{eq:single-mode}
u_\infty^-(\theta)=e^{im\theta},\qquad m\in\mathbb Z.
\end{equation}
We consider scattering at the unit disk. This simplified setup has explicit solutions available due to the rotational invariance of the problem. As a consequence of the unit amplitude for the incident field, the scattering amplitudes are simply the far-field amplitudes of the scattered fields. Single mode data also have the advantage of isolating the centrifugal barrier as \(|m|\) approaches \(k\). We discuss numerical and analytical solutions for constant (Sec.~\ref{sec:single-mode-constant}), short-range (Sec.~\ref{sec:single-mode-short-range}), and long-range (Sec.~\ref{sec:single-mode-long-range}) media.

\subsubsection{Constant media}\label{sec:single-mode-constant}

\paragraph{Explicit solution}\label{sec:explicit-constant}
Seeking modal solutions \(U(r,\theta)=U_m(r)\,e^{im\theta}\) with \(m\in\mathbb{Z}\),
we obtain the radial equation
\begin{equation}\label{eq:modal-constant}
U_m'' + \frac{1}{r}U_m' + \left(k^2 - \frac{m^2}{r^2}\right)U_m = 0.
\end{equation}
The solution can be written as a linear combination of Hankel functions of order $m$:
\[
U_m(r) = U_m^-(r)+U_m^+(r) = A_m\,H_{m}^{(2)}(kr)+B_m\,H_{m}^{(1)}(kr).
\]
The coefficients $(A_m,B_m)$ are determined by the boundary conditions. We consider single-mode incidence from infinity, which determines $A_m$, and a sound-soft unit disk, which determines $B_m$. To express the modal scattering amplitude, we write the asymptotic behavior of the Hankel functions for large argument \cite{DLMF, olver1997asymptotics}:
\begin{equation}\label{eq:hankel-asymptotics}
H_m^{(2)}(kr) \sim \frac{e^{-ikr}}{\sqrt{r}} \sqrt{\frac{2}{\pi k}} e^{i \frac{\pi}{2} \left(m+\frac{1}{2}\right)},  \qquad 
H_m^{(1)}(kr) \sim \frac{e^{ikr}}{\sqrt{r}} \sqrt{\frac{2}{\pi k}} e^{-i \frac{\pi}{2} \left(m+\frac{1}{2}\right)}.
\end{equation}
The condition \eqref{eq:single-mode} together with the definition of the far-field amplitude \eqref{eq:radiation-expansion} implies
\begin{equation}\label{eq:Am}
A_m = \sqrt{\frac{\pi k}{2}}\,e^{-i \frac{\pi}{2} \left(m+\frac{1}{2}\right)},
\end{equation}
so that \(U_m^-(r,\theta)\sim r^{-1/2}e^{-ikr}u_\infty^-(\theta)\) as \(r\to\infty\). Imposing the sound-soft condition at a disk of unit radius $r=1$, we have $U_m(1)=U_m^-(1)+U_m^+(1)=0$, so
\[
A_m\,H_m^{(2)}(k)+B_m\,H_m^{(1)}(k)=0
\quad\Longrightarrow\quad
B_m=-A_m\,\frac{H_m^{(2)}(k)}{H_m^{(1)}(k)}.
\]
The scattered solution for a sound-soft disk of unit radius is
\begin{equation}\label{eq:scattered-general}
U_m^+(r) =- A_m \frac{H_m^{(2)}(k)}{H_m^{(1)}(k)}\,H_m^{(1)}(kr), \qquad r>1.
\end{equation}

The scattering map at wavenumber \(k\) is the operator $S(k): u_\infty^-(\theta) \mapsto u_\infty^+(\theta)$, mapping incoming radiation to outgoing radiation. Due to the rotational invariance of our problem, the scattering map is diagonal. In the single-mode case, the modal coefficients are given by the far-field pattern of the scattered solution, 
which follows from the asymptotic behavior of Hankel functions in \eqref{eq:hankel-asymptotics}:
\begin{equation}\label{eq:S-constant}
S_m(k) = i(-1)^m\,\frac{H_m^{(2)}(k)}{H_m^{(1)}(k)}.
\end{equation}
The scattering map is unitary, \(|S_m|=1\).

\paragraph{Modal numerical solution and the centrifugal barrier}
We solve the modal constant-index Helmholtz equation by expanding the NIC solution into Fourier modes \(u^\pm(\rho,\theta)=\sum_{m\in\mathbb{Z}} u_m^\pm(\rho)e^{im\theta}\). Denoting the radial derivative by a prime, we get from \eqref{eq:compactified} the radial ordinary differential equation (ODE),
\begin{equation}\label{eq:mode-ode}
\left((1-\rho)^2 (u^\pm_m)'\right)' \pm 2ik\,(u^\pm_m)' - \frac{m^2-\frac14}{\rho^2}\,u^\pm_m = 0, 
\qquad \rho\in(\rho_\Gamma,\rho_{\scri})=\left(\frac12, 1\right).
\end{equation}
The minus sign corresponds to the incoming branch \(u_m^-\) and the
plus sign to the outgoing branch \(u_m^+\). 
The purely incoming physical solution is singular at the origin. 
We solve the above equation on the punctured domain \(\rho\in(\rho_\Gamma,1)\), which corresponds to the physical exterior \(r\in(1,\infty)\), using a standard 1D Chebyshev collocation method with Gauss--Lobatto points.

For code verification and to investigate the impact of the centrifugal barrier separately on the incident and scattered solutions, we impose the analytic inner boundary data for the scattered field directly from the analytic solution\footnote{This decoupling (analytic inner boundary for \(u^+\)) is used only in the 1D case to test the accuracy of incident and scattered solvers independently. In the 2D experiments, we couple the two solutions at the obstacle boundary numerically via the sound-soft condition using \eqref{eq:numerical_boundary_data} or \eqref{eq:boundary-variable}.} given in \eqref{eq:scattered-general}. The boundary data for the incident and scattered fields are
\[
u^-_m|_{\scri^-}=1, \qquad 
u^+_m|_\Gamma 
=-e^{-ik}A_m H_m^{(2)}(k),
\]
with $A_m$ given in \eqref{eq:Am}.

\begin{figure}[t]
  \centering
  \includegraphics[width=0.9\textwidth]{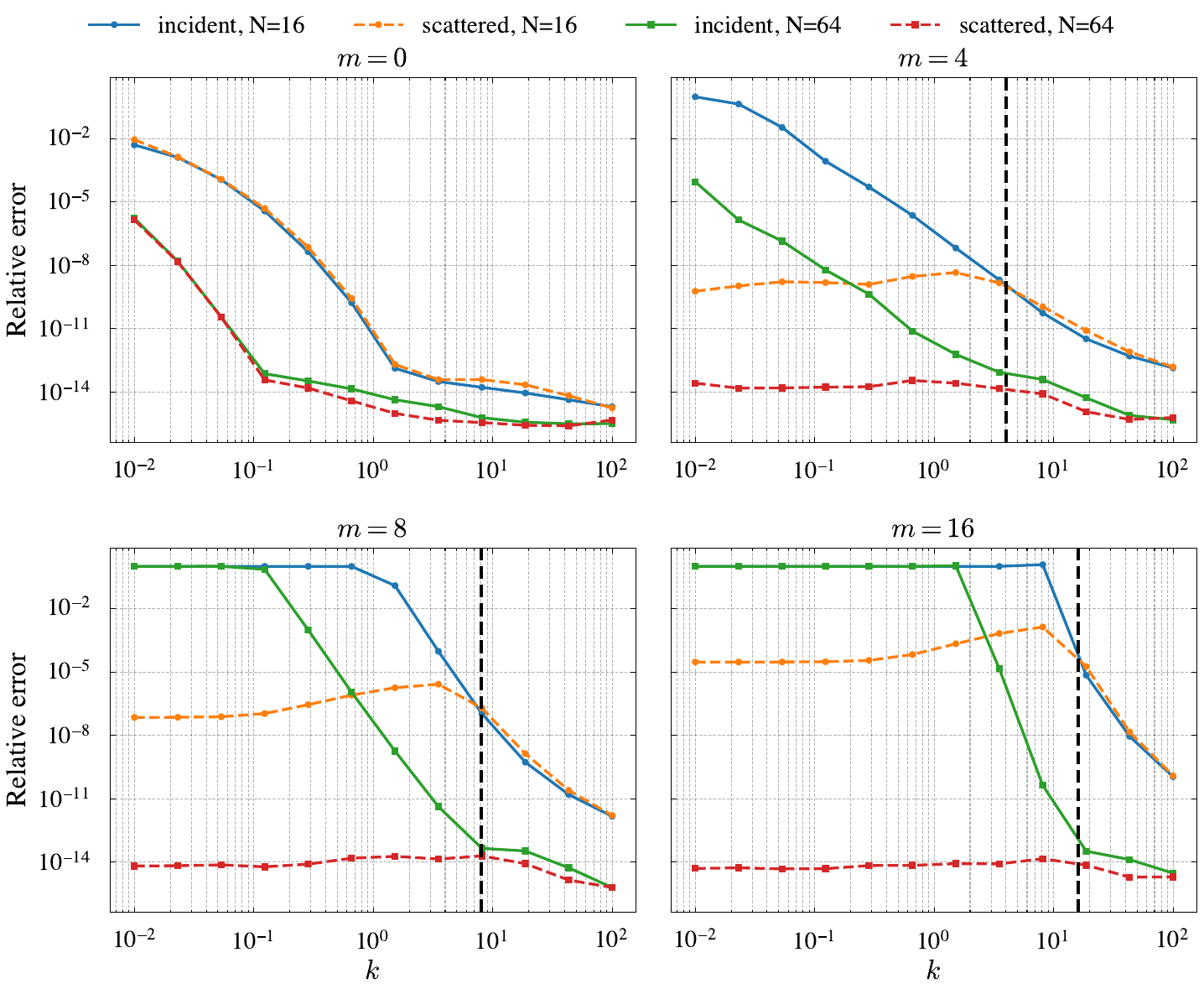}
  \caption{Relative errors in the discrete $\ell^2$ norm of incident (solid) and scattered (dashed) fields vs.\ wavenumber \(k\) for four azimuthal indices \(m\in\{0,4,8,16\}\) in two resolutions $N=\{16,64\}$ using a 1D Chebyshev spectral method. The vertical dashed lines mark the  centrifugal-barrier. The scattered field is systematically more accurate; incident errors are high below the barrier but drop rapidly beyond.}
  \label{fig:km-errors}
\end{figure}

In Fig.~\ref{fig:km-errors} we show the relative errors of the incident and scattered fields as a function of wavenumber \(k\in[10^{-3},10^3]\) for the azimuthal indices \(m=\{0,4,8,16\}\). 
The relative error between the numerical solution $u^\text{num}$ and the explicit solution $u$ is computed as the discrete $\ell^2$ norm on the Gauss--Lobatto grid:
\begin{equation}
\label{eq:discrete-l2-error}
\mathcal{E}_{\mathrm{rel}}
= \frac{\|u^\text{num} - u\|_{\ell^2}}{\|u\|_{\ell^2}}, 
\quad \text{where} \quad 
\|u\|_{\ell^2}^2 := \sum_{j=0}^N |u(\rho_j)|^2.
\end{equation}
The most prominent feature in the accuracy of modal solutions is the centrifugal barrier. For the \(m\)-mode, the physical radial wavenumber is \(k^2-m^2/r^2\). For large \(m\), the effective potential \(m^2/r^2\) limits waves toward small \(r\) creating a centrifugal barrier near \(m/k\), marked by vertical dashed lines in Fig.~\ref{fig:km-errors}. The numerical solution demonstrates spectral accuracy away from the centrifugal barrier but not below. Below the barrier (\(r<m/k\)), incoming single modes are ill-conditioned because the purely incoming branch \(H_m^{(2)}(z)\) grows as $z^{-m}$ for \(z=kr\to0\). 
In contrast, the scattered problem is well-conditioned for small $k$ because the scattering map \eqref{eq:S-constant} behaves as $S_m \sim i(-1)^{m+1} \left(1+ \mathcal{O}\left(k^{2m}\right)\right)$
so the scattered modal amplitude decays rapidly in \(m\) at small \(k\). Therefore, the scattered field is consistently more accurate than the incident field for $m\ne 0$. 

\paragraph{Full numerical implementation and convergence}
We solve the NIC-Helmholtz equation \eqref{eq:compactified} with single-mode incident boundary data \eqref{eq:single-mode} using a spectral method implemented in the Dedalus framework \cite{burns2020dedalus}. The solution is represented in a Fourier--Chebyshev basis,
\vspace{-3mm}
\[
u^\pm(\rho, \theta) \approx \sum_{n=0}^{N} \sum_{m=-M}^{M} \hat{u}^\pm_{n,m} T_n(\rho) e^{i m \theta},
\vspace{-3mm}
\]
where $ T_n(\rho) $ are Chebyshev polynomials and $ e^{i m \theta} $ are Fourier modes. In Dedalus, the Dirichlet boundary enforced via a tau formulation in the radial direction, and the solution is evaluated on Chebyshev grids that include the interval endpoints.

\begin{figure}[t]
  \centering
  \includegraphics[width=0.48\textwidth]{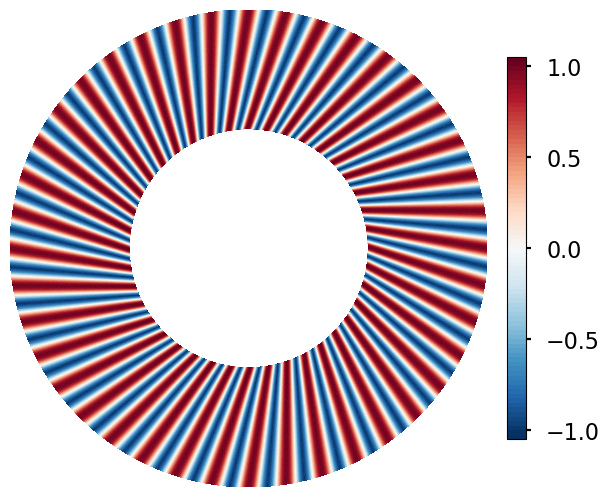}\hfill
  \includegraphics[width=0.48\textwidth]{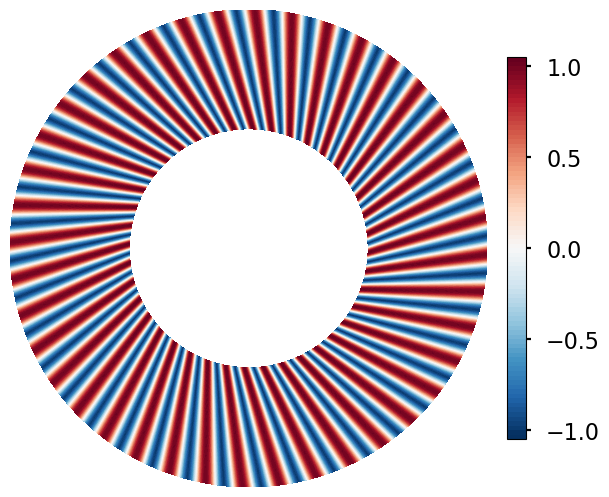}
  \caption{The rescaled incident and scattered fields, $u^\pm$, for single-mode incoming radiation with $m=40$ and $k=120$. The inner circle is the obstacle surface; the outer circle is null infinity (compare Fig.~\ref{fig:incident-scattered}). The solution is dominated by the angular oscillation with relatively straight rays connecting the obstacle with infinity. The scattering map is unitary (pure phase).}
  \label{fig:modal-solutions}
\end{figure}

Figure~\ref{fig:modal-solutions} shows the real part of the rescaled solutions \(u^\pm(\rho,\theta)\) for the incident and scattered fields with \(m=40\) and \(k=120\). 
The inner circle represents the obstacle surface \(\rho_\Gamma=\frac12\) and the outer circle represents null infinity \(\rho_{\scri^\pm}=1\). We extract the trace of the incident field on the obstacle from the incoming Helmholtz problem, which provides the numerical Dirichlet boundary data for the scattered field through \eqref{eq:numerical_boundary_data}. 

The solution is dominated by the angular oscillation \(e^{i m\theta}\) (forty oscillations around the circle). The variation in the compactified radial direction is mild. The bands form almost straight rays connecting \(\scri^\pm\) to the obstacle, indicating that NIC transports information essentially along straight characteristics across the exterior domain. The only visible deviation from straight rays is a slight bending of the phase lines, caused by the centrifugal barrier for \(m=40\). Both the incident and scattered rescaled fields remain smooth and of order one from the obstacle to null infinity, despite the relatively high wavenumber \(k=120\). 

\begin{figure}[t]
  \centering
  \includegraphics[width=0.95\textwidth]{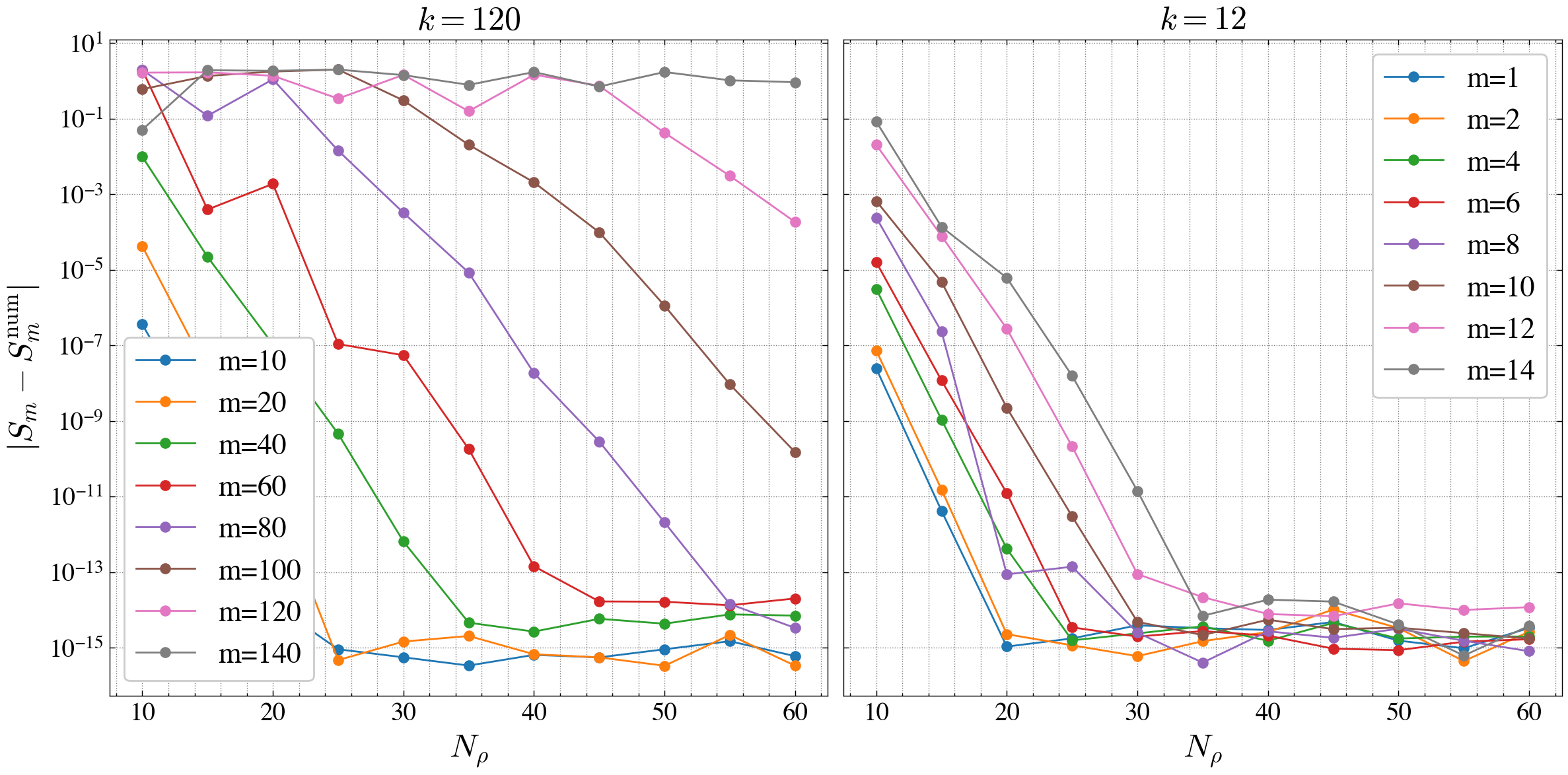}
\caption{Radial spectral convergence of the modal scattering map at null infinity. 
We plot on a semilog scale the absolute error \(|S_m(k)-S_m^{\mathrm{num}}(k)|\) against the radial polynomial degree \(N_\rho\) for two wavenumbers (left: \(k=120\); right: \(k=12\)) using the same set of ratios \(m/k\), and fixed angular resolution \(N_\theta=2m+5\). 
For all cases the error decays exponentially with \(N_\rho\) until it reaches saturation. The saturation level and the total error increase with \(m/k\) and $k$.}
  \label{fig:convergence}
\end{figure}

We now compute the modal scattering map at null infinity for single modes to check the numerical convergence of the coupled solvers. Figure \ref{fig:convergence} compares the numerical value with the explicit expression given in \eqref{eq:S-constant}. The error in the scattering map is defined as
\[
\mathcal{E}_m(k)=\bigl|\,S_m(k)-S_m^{\text{num}}(k)\,\bigr|.
\]
The absolute difference measures the phase error because the exact scattering coefficient satisfies $|S_m(k)|=1$. As both incident and scattered solvers use spectral resolution, we expect exponential decay of the error with increasing resolution, as confirmed by Fig.~\ref{fig:convergence}. In the experiments, the angular resolution is fixed to \(N_\theta=2m+5\). The two panels use the same set of ratios \(m/k\) (left: \(k=120\), right: \(k=12\)), so each curve represents the same centrifugal-barrier. In both panels the error decays exponentially in \(N_\rho\) until a saturation level that increases with \(m/k\). As \(m/k\) grows and the turning point reaches the obstacle, the onset of convergence is delayed. At fixed \(m/k\) the errors are larger for \(k=120\) than for \(k=12\), consistent with the sharpening of the transition layer near the turning point \cite{olver1997asymptotics}. In summary, we observe the expected spectral convergence in radial resolution. In accordance with standard WKB analysis, the error increases with $m$ for fixed $k$, and with $k$ for fixed $m/k$.

\subsubsection{Short-range media} \label{sec:single-mode-short-range}
We compute the scattering map for a sound-soft unit disk embedded in a short-range medium. The perturbation of the refractive index is radial with unbounded support and quadratic decay:
\begin{equation}\label{eqn:quadratic}
n^2(r) = 1 + \frac{\kappa^2}{r^2}.
\end{equation}
\begin{figure}[t]
  \centering
  \includegraphics[width=0.95\textwidth]{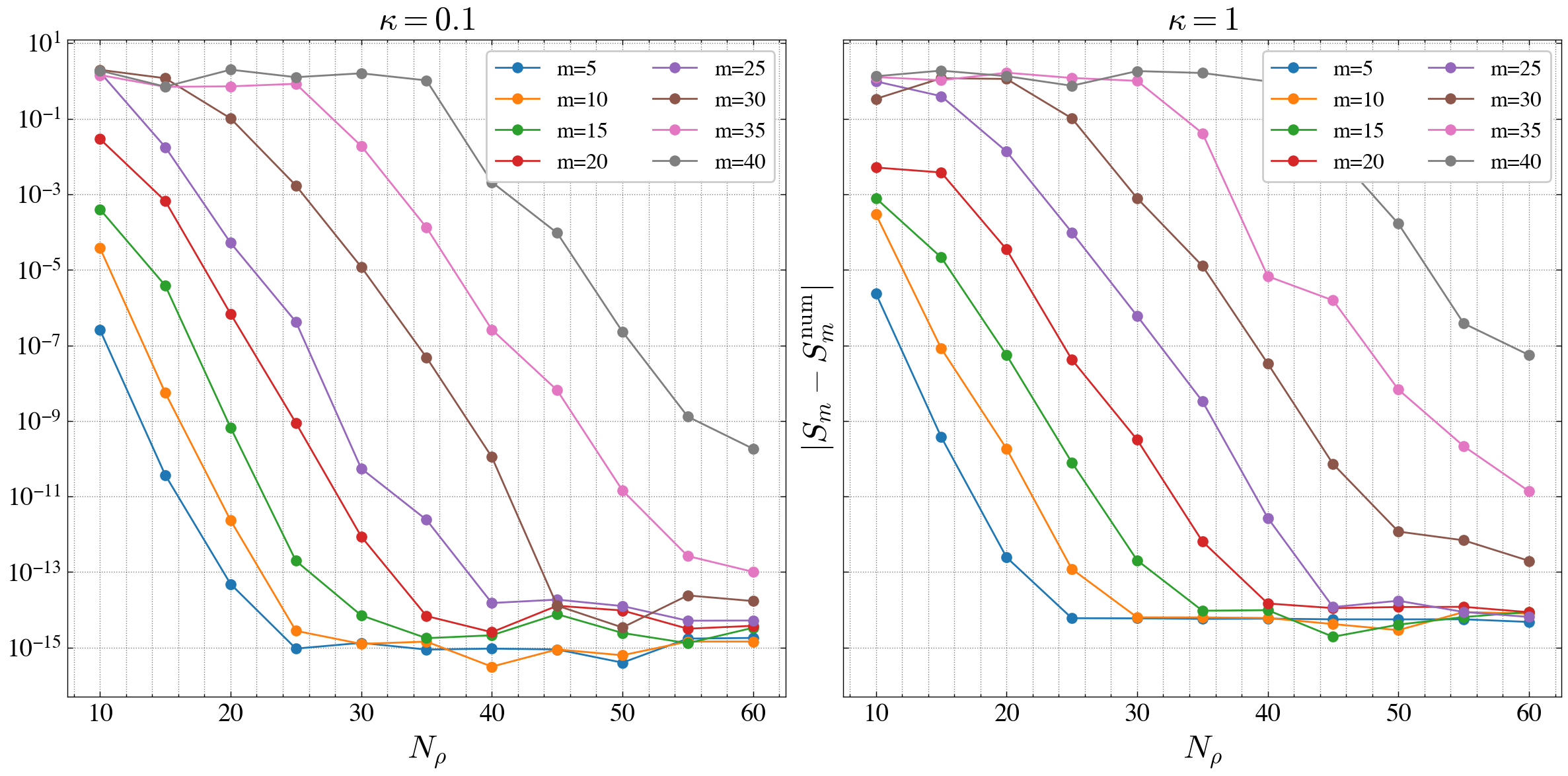}
\caption{Radial spectral convergence of the modal scattering map at null infinity for variable media with a radially symmetric refractive index of unbounded support and quadratic decay \eqref{eqn:quadratic}.
We plot on a semilog scale the absolute error \(|S_m(k)-S_m^{\mathrm{num}}(k)|\) against the radial polynomial degree \(N_\rho\) for two values of the strength parameter (left: \(\kappa=0.01\); right: \(\kappa=0.99\)) using the same set of $m$-modes, same wave number $k=20$, and angular resolution \(N_\theta=2m+5\). The error decays exponentially with \(N_\rho\) until it reaches saturation. 
}
  \label{fig:convergence_quadratic}
\end{figure}
We drive the problem by a single incoming \(m\)-mode prescribed at infinity. We can explicitly solve the modal equation 
\[
U_m'' + \frac{1}{r}U_m' + \left(k^2 - \frac{\nu^2}{r^2}\right)U_m = 0,
\qquad
\nu^2 := m^2 - \kappa^2 k^2.
\]
The solution can be written as a linear combination of Hankel functions of order $\nu$: 
\[
U_m(r) = A_\nu\,H_{\nu}^{(2)}(kr)+B_\nu\,H_{\nu}^{(1)}(kr).
\]
Imposing the incoming single-mode radiation condition \eqref{eq:single-mode} and proceeding as in Sec.~\ref{sec:explicit-constant}, we obtain the modal scattering map:
\begin{equation}\label{eq:Sm-quadratic}
S_m(k;\kappa) = i(-1)^\nu\,\frac{H_{\nu}^{(2)}(k)}{H_{\nu}^{(1)}(k)},
\qquad
\nu = \sqrt{\,m^2-\kappa^2 k^2\,}\,.
\end{equation}

For the numerical solver, we use characteristic coordinates and the rescaling \eqref{eq:u-variable-rescaled}. The resulting compactified equation \eqref{eq:compactified-variable} has regular coefficients at null infinity and is solved numerically using the same spectral--Galerkin method as in the previous section. The rescaled fields are coupled through the boundary condition \eqref{eq:boundary-variable} where
\begin{equation} \label{eq:integral-quadratic}
\mathcal{N}(r) = \int\sqrt{1+\frac{\kappa^2}{r^2}}\, dr = \sqrt{r^2+\kappa^2} - \kappa \operatorname{arcsinh} \frac{\kappa}{r}.
\end{equation}
Spectral convergence of the numerical scattering map to the explicit solution \eqref{eq:Sm-quadratic} is demonstrated in Fig.~\ref{fig:convergence_quadratic} for two separate values of the strength parameter \(\kappa\). We observe exponential decay of the error with increasing radial resolution \(N_\rho\) until saturation. The numerical errors are mildly higher for a larger strength parameter.

It is notable that NIC does not rely on the availability of the fundamental solution and can therefore be applied to variable media with minimal changes to the formalism. The main ingredient when using characteristic coordinates is the solution of the eikonal equation as in \eqref{eq:integral-quadratic}. However, NIC only requires a regular compactification at infinity and therefore does not depend on the explicit knowledge of the eikonal solution at every point in the exterior domain. In Sec.~\ref{sec:hyperbolic}, we provide the computation of a variable media solution in hyperbolic coordinates as an example for a setup that does not require characteristic coordinates. 

\subsubsection{Long-range media}\label{sec:single-mode-long-range}
A stronger modification of the method is required for long-range media, where the refractive index has a Coulomb-type tail. We consider a radial refractive index with perturbations of linear decay
\begin{equation}\label{eqn:linear}
n^2(r)=1+\frac{\kappa}{r}.
\end{equation}
To construct the explicit solution analytically, we solve the radial ODE
\begin{equation}\label{eq:modal-long-range}
U_m''(r) + \frac{1}{r}U_m'(r)
+ \left(k^2 + \frac{k^2\kappa}{r} - \frac{m^2}{r^2}\right)U_m(r)=0.
\end{equation}
Related solutions for radially varying acoustic media can be found in \cite{martin2002acoustic}. Introducing the scaled radial variable $z=kr$ and applying a Liouville transform via $U_m(r) = z^{-1/2} w(z)$, we obtain the  Coulomb radial equation
\[
w''(z) + \left(1 - \frac{2\eta}{z} - \frac{\lambda(\lambda+1)}{z^2}\right)w(z)=0,
\]
with parameters
\begin{equation}
\lambda = m-\frac12,\qquad \lambda(\lambda+1)=m^2-\frac14,
\qquad
\eta = -\frac{k\kappa}{2}.
\end{equation}
The solution can be written in terms of the Coulomb-Hankel functions \cite{DLMF} as
\begin{equation}\label{eq:general-mode}
U_m(r)
= \frac{1}{\sqrt{k r}}\left(
A_\lambda\,H^{(-)}_{\lambda}(\eta,k r)
+ B_\lambda\,H^{(+)}_{\lambda}(\eta,k r)
\right).
\end{equation}
To determine the coefficients, we consider the asymptotic behavior of the Coulomb-Hankel functions as $z\to\infty$. Defining the Coulomb phase shift $\sigma_\lambda(\eta)=\arg\Gamma(\lambda+1+i\eta)$, we have the following asymptotics for large argument \cite{DLMF}:
\begin{align}\label{eq:LR-coulomb-asymp}
H^{(\pm)}_{\lambda}(\eta,z)
&\sim \exp\!\left( \pm i\left[z-\eta\log(2z)-\frac{\pi}{2}\lambda+\sigma_\lambda(\eta)\right] \right), \nonumber \\
&\sim \exp\!\left(\pm ik \left[r + \frac{\kappa}{2}\log r\right] \right)\,
  \exp\!\left(\mp i\left[\frac{k\kappa}{2}\log(2k) +\frac{\pi}{2}\lambda-\sigma_\lambda(\eta)\right] \right).
\end{align}
We define the long-range phase
\[
\Phi(r;\kappa):=r+\frac{\kappa}{2}\log r,
\]
so that $e^{\pm i k \Phi}$ captures the $r$-dependent part of the Coulomb phase. The modal coefficients have the long-range asymptotic expansion
\begin{equation}
U_m(r)\sim
\frac{1}{\sqrt{r}}\left(
e^{-ik \Phi(r;\kappa)}\,u_{m,\infty}^{-} + e^{+i k\Phi(r;\kappa)}\,u_{m,\infty}^{+} \right).
\end{equation}
We prescribe single-mode incidence with unit amplitude at infinity as in \eqref{eq:single-mode}. This condition implies that the far-field amplitude of the scattered field is the scattering amplitude:
\begin{equation}\label{eq:LR-far-field}
u_{m,\infty}^{-}=1 \ \implies \
u_{m,\infty}^{+}=S_m(k;\kappa)
\end{equation}
To normalize the incident solution such that incoming radiation data has unit amplitude at infinity, we set
\begin{equation}\label{eq:LR-A-lambda}
A_\lambda :=
\sqrt{k}\, 
\exp\!\left(
-i\left[\eta\log(2k)+\frac{\pi}{2}\lambda-\sigma_\lambda(\eta)\right]
\right).
\end{equation}
With this choice the incident field
\begin{equation}\label{eq:LR-Uinc}
U^-_m(r)
:=A_\lambda\,
\frac{H^{(-)}_{\lambda}(\eta,kr)}{\sqrt{kr}}
\end{equation}
satisfies
\begin{equation}
\lim_{r\to\infty}
\sqrt{r}\,e^{+i k \Phi(r;\kappa)}\,
U^{\mathrm{inc}}_m(r,\theta)=u_{m,\infty}^{-} = 1.
\end{equation}
For sound-soft scattering, we impose the homogeneous Dirichlet boundary condition on the total field at the unit circle: $U^-_m(1)=U^-_m(1)+U^+_m(1)$. The condition determines the scattered mode
\begin{equation}\label{eq:LR-Uscat}
U^+_m(r)
:=
-\,A_\lambda\,
\frac{H^{(-)}_{\lambda}(\eta,k)}{H^{(+)}_{\lambda}(\eta,k)}\,
\frac{H^{(+)}_{\lambda}(\eta,kr)}{\sqrt{kr}}.
\end{equation}
The Coulomb-modified scattering amplitudes are then obtained by taking the outgoing far-field limit:
\begin{align}
S_m(k;\kappa)
&= \lim_{r\to\infty} \sqrt{r}\,e^{-i k \Phi(r;\kappa)}\, U^+_m(r),\\
&=-\,\frac{A_{\lambda}}{\sqrt{k}}\,
\frac{H^{(-)}_{\lambda}(\eta,k)}{H^{(+)}_{\lambda}(\eta,k)}\,
\exp\!\left(
i\Big[\frac{k\kappa}{2}\log(2k)-\frac{\pi}{2}\lambda+\sigma_{\lambda}(\eta)\Big]\right),
\end{align}
which, using \eqref{eq:LR-coulomb-asymp} and \eqref{eq:LR-A-lambda}, yields the explicit expression
\begin{equation}\label{eq:Sm-long-range}
S_{m}(k;\kappa)
=
- i (-1)^m\,\frac{H^{(-)}_{\lambda}(\eta,k)}{H^{(+)}_{\lambda}(\eta,k)}\,
e^{i k\kappa\log(2k)+ 2 i \sigma_{\lambda}(\eta)}\,.
\end{equation}
Comparing this expression with the constant-index \eqref{eq:S-constant} and the short-range \eqref{eq:Sm-quadratic} amplitudes, we see that the long-range scattering amplitude is modified by a Coulomb phase factor that depends on the strength parameter \(\kappa\). This relatively simple formula demonstrates that, for long-range media, an accurate scattering computation necessarily involves both the Coulomb-type background phase and the obstacle response. This calculation is analytically cumbersome with plane wave incidence. The standard treatment described in Sec.~\ref{sec:standard_variable}, typically assuming plane wave incidence, results in the variable-coefficient Helmholtz equation \eqref{eq:sourced-helmholtz}. Modal solutions are given by nontrivial integral representations involving Coulomb-Hankel functions and Bessel functions with no simple closed form for the scattering map. By contrast, prescribing single-mode incident data at infinity leads to the homogeneous Coulomb equation \eqref{eq:modal-long-range} and the scattering map reduces to a single explicit phase factor.

\begin{figure}[t]
  \centering
  \includegraphics[width=0.95\textwidth]{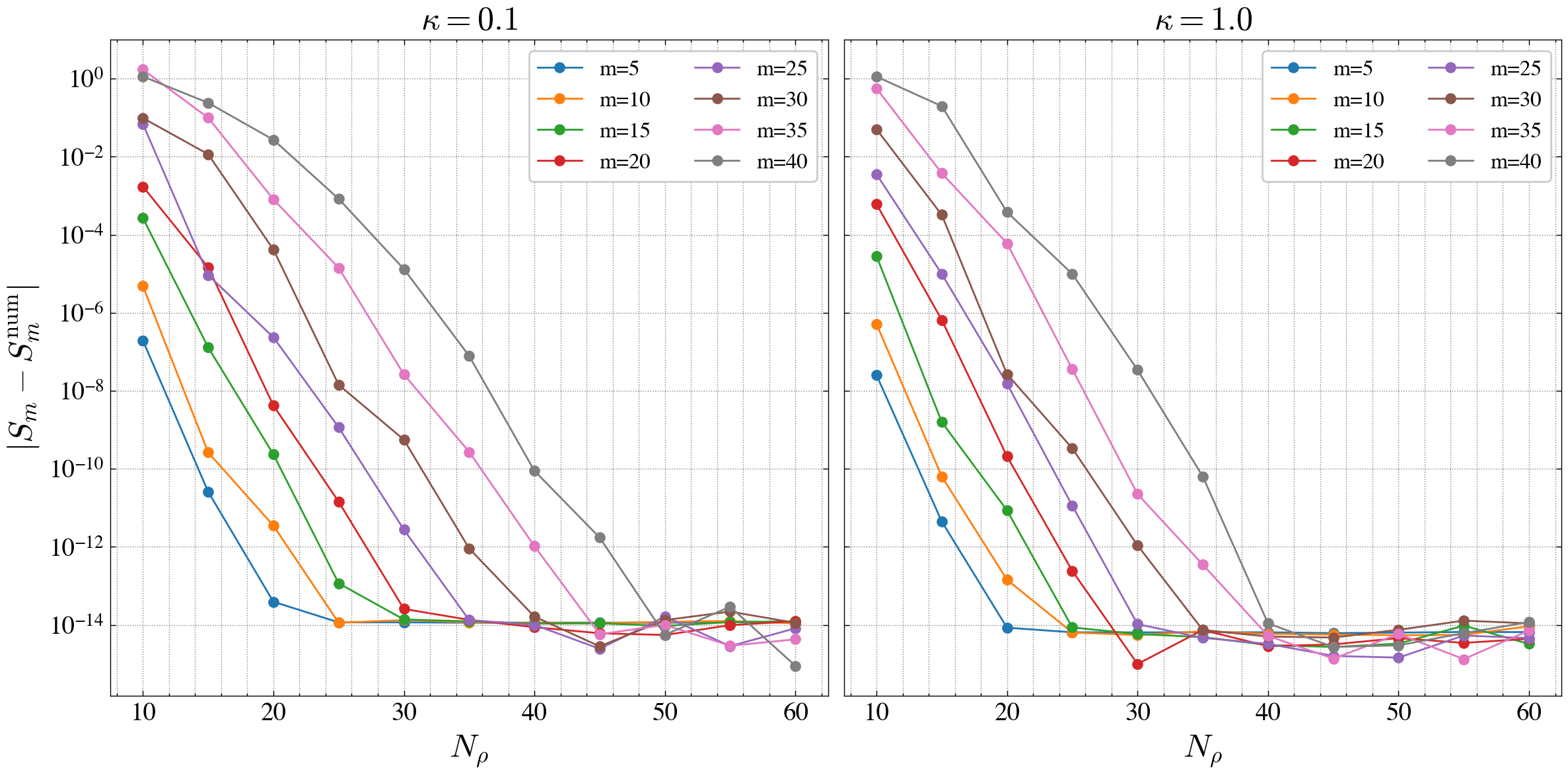}
\caption{Radial spectral convergence of the modal scattering map at null infinity for variable media with a radially symmetric refractive index of unbounded support whose deviation from unity decays quadratically \eqref{eqn:quadratic}.
We plot on a semilog scale the absolute error \(|S_m(k)-S_m^{\mathrm{num}}(k)|\) against the radial polynomial degree \(N_\rho\) for two values of the strength parameter (left: \(\kappa=0.1\); right: \(\kappa=0.9\)) using the same set of $m$-modes, same wavenumber $k=40$, and angular resolution \(N_\theta=2m+5\). For all cases the error decays exponentially with \(N_\rho\) until it reaches saturation. The saturation level and the total error are mildly higher for a larger strength parameter.}
  \label{fig:convergence_linear}
\end{figure}

For the NIC implementation, we solve the eikonal equation,
\begin{equation}
\mathcal{N}(r)
=\int\sqrt{1+\frac{\kappa}{r}}\,dr 
= \sqrt{r(r+\kappa)} - \kappa \log\!\left(\sqrt{r+\kappa}-\sqrt{r}\right) - C(\kappa),
\end{equation}
where the integration constant,
$C(\kappa) := \frac{1}{2}\left(\kappa + \ln \frac{2}{\kappa} \right)$,
is chosen so that as $r\to\infty$, we have
$\mathcal{N}(r) \sim r + \frac{\kappa}{2}\log r = \Phi(r;\kappa)$. The phase and the normalization of the solution ensure that the incident NIC field, 
$u^-(\rho,\theta) := \sqrt{r(\rho)}\,
\exp\!\left(i k\,\mathcal{N}(r(\rho))\right)\,
U^-_m(r(\rho),\theta)$ satisfies $u^-(1,\theta)=e^{i m \theta}$ at null infinity and the scattered NIC field has the far-field corresponding to the scattering amplitudes \eqref{eq:Sm-long-range}.

Spectral convergence of the numerical scattering map to the explicit solution \eqref{eq:Sm-long-range} is demonstrated in Fig.~\ref{fig:convergence_linear} for two separate values of the strength parameter \(\kappa\). As in Fig.~\ref{fig:convergence_quadratic}, we observe exponential decay of the error with increasing radial resolution \(N_\rho\) until saturation. The numerical errors are mildly higher for a larger strength parameter.


\subsection{Gaussian beam}\label{sec:gaussian}
The single-mode incoming radiation discussed in the previous section is mathematically appealing, but does not represent a field radiated by a localized source. To localize the incoming radiation at infinity in the angular direction, we prescribe the data, \(u_\infty^-(\theta)\), as a narrow Gaussian. Such data generates a beam-like incident field.

A circular Gaussian beam, plotted on the right panel of Fig.~\ref{fig:incoming-panel}, provides smooth angular data with exponential modal decay,
\begin{equation}\label{eq:gaussian}
u_\infty^-(\theta)=\frac{1}{\sqrt{2\pi I_0(2\beta)}}\exp\left(\beta\cos(\theta-\theta_0)\right),
\end{equation}
where \(\theta_0\) is the beam direction and \(I_0\) is the modified Bessel function of the first kind ensuring that the incident field is normalized, \(\|u_\infty^-\|_{L^2(0,2\pi)}=1\) so that the total energy flux of the incident field is unity.
The parameter \(\beta>0\) is the concentration, governing the angular spread of the beam: small~\(\beta\) yields an almost isotropic incidence, while large~\(\beta\) produces a narrowly focused beam around~\(\theta_0\). The square modulus \(|u_\infty^-(\theta)|^2\) is essentially the von Mises distribution, a circular analogue of the normal distribution \cite{watson1982distributions}. The beam has exponentially decaying modal coefficients, which reduces Gibbs oscillations and improves numerical stability. Using the Jacobi-Anger expansion \(e^{\beta\cos\theta}=\sum_{n\in\mathbb{Z}} I_n(\beta)\,e^{in\theta}\), we compute the modal spectrum as
\begin{equation}\label{eq:VM-am}
u^-_{m,\infty}=\frac{1}{\sqrt{2\pi I_0(2\beta)}}\,I_{m}(\beta)\,e^{-im\theta_0}.
\end{equation}
The modes decay exponentially for large mode number as $|m|\to\infty$ with the estimate $|u^-_{m,\infty}| \sim e^{-m^2/(2\beta)}$  \cite{olver1997asymptotics}. The Gaussian beam is concentrated for larger values of $\beta$ exciting higher $m$-modes.


\subsubsection{Constant media}

Given the modal representation \eqref{eq:VM-am} and the single-mode solution for the sound-soft unit disk discussed in Sec.~\ref{sec:explicit-constant}, we can construct an explicit solution for the scattered field. Each mode is multiplied by the corresponding  scattering amplitude: $u^+_{m,\infty} = S_m(k)\,u^-_{m,\infty}$, where the amplitude for the constant-index case is given in \eqref{eq:S-constant}. The exact scattered far-field modal amplitudes read
\begin{equation}\label{eq:VM-scattered-modes}
u^+_{m,\infty}
= S_m(k)\,u^-_{m,\infty}
= i(-1)^m\,\frac{H_m^{(2)}(k)}{H_m^{(1)}(k)}\,
\frac{I_m(\beta)}{\sqrt{2\pi I_0(2\beta)}}\,e^{-im\theta_0}.
\end{equation}
We obtain the corresponding scattered far-field pattern by summing up the Fourier series, $u_\infty^+(\theta)
= \sum_{m\in\mathbb{Z}} u^+_{m,\infty} \, e^{im\theta}$.
This explicit representation provides a rigorous benchmark for the numerical scheme. In the implementation we work in the Fourier space for the angular discretization, with integer mode indices $m \in \{-N_\theta/2,\dots,N_\theta/2-1\}$. For a given $(\beta,k)$ we first form the exact incoming modal amplitudes $u^-_{m,\infty}$ using \eqref{eq:VM-am} and then multiply by $S_m(k)$ to obtain the scattered modes $u^+_{m,\infty}$ in \eqref{eq:VM-scattered-modes}. Assigning these coefficients to a spectral vector and applying a backward transform yields the exact scattered far field $u_{\infty}^+(\theta)$ on the same discrete angular grid as the numerical solution. The numerical solver produces discrete incident and scattered solutions $u^-$ and $u^+$ from \eqref{eq:compactified}. We evaluate these fields at null infinity, which corresponds to the outer boundary of the domain, $\rho=1$, sampled on the discrete grid $\{\theta_j\}_{j=0}^{N_\theta-1}$.

\begin{figure}[t]
  \centering
  \includegraphics[width=0.49\textwidth]{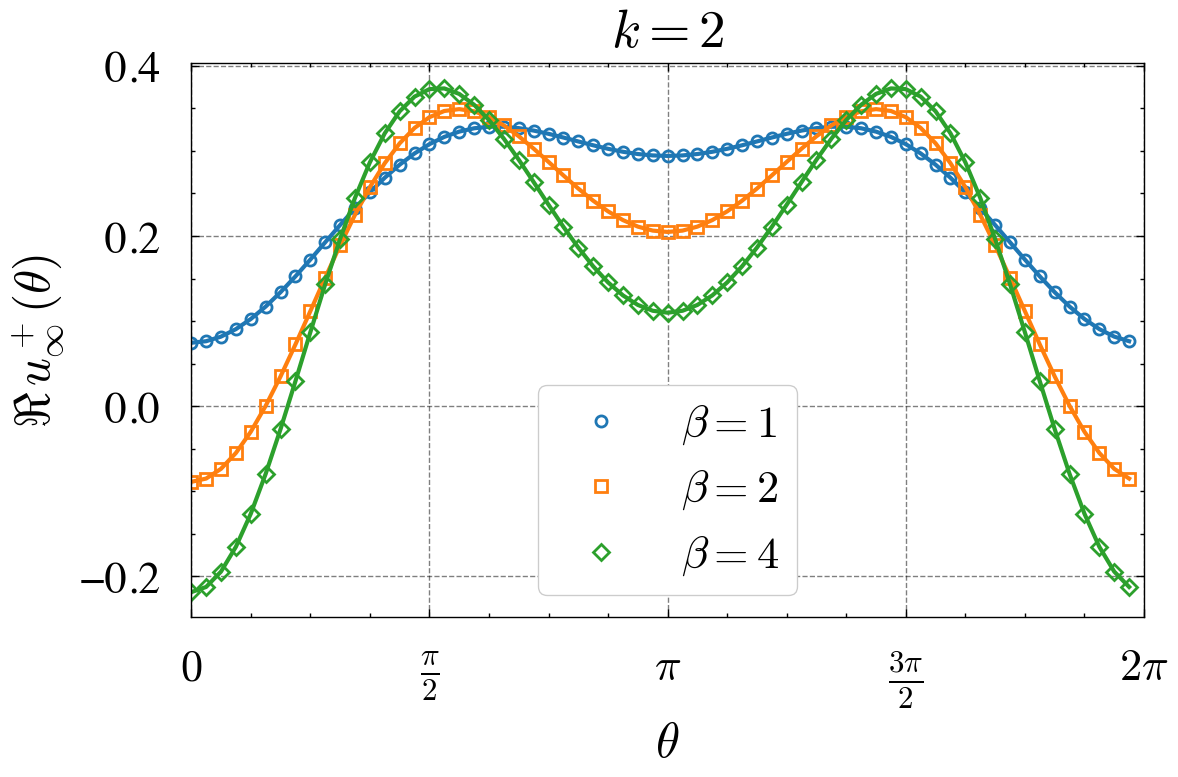}\hfill
  \includegraphics[width=0.49\textwidth]{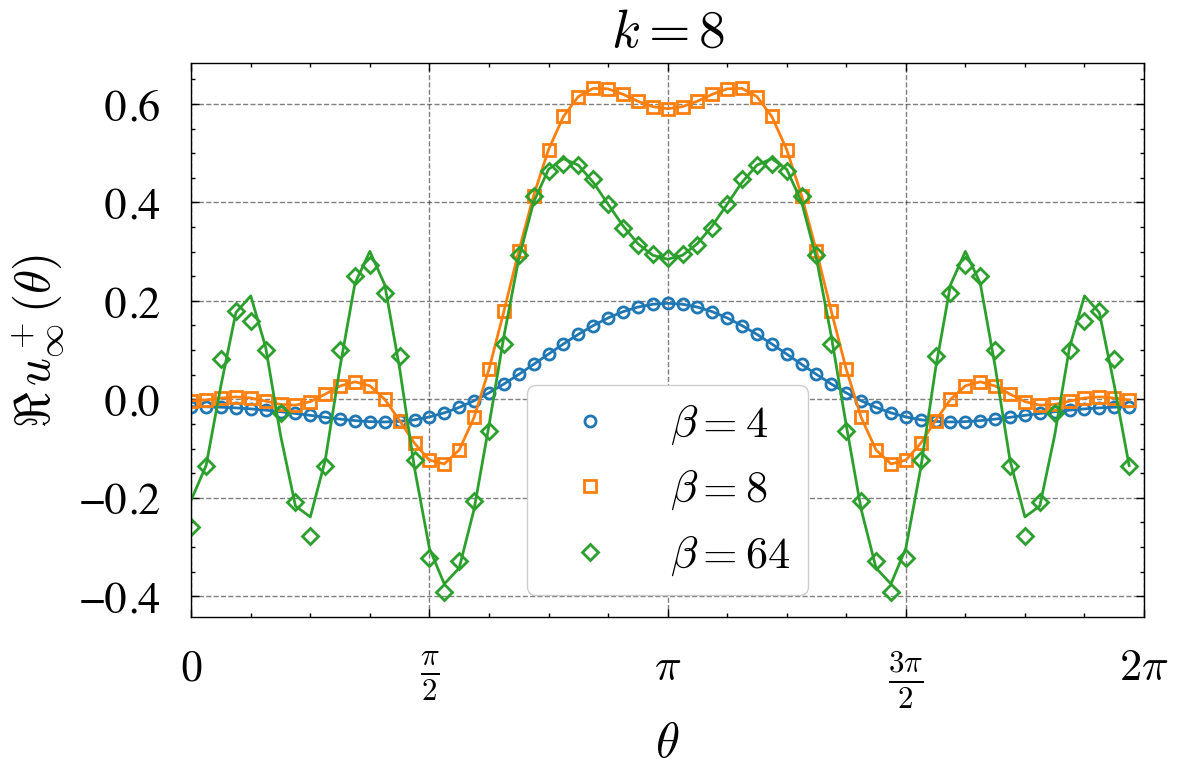}
\caption{Comparison of the exact scattered far field (solid lines) and the numerical solution (open markers) for Gaussian beam incidence \eqref{eq:gaussian} at different concentrations $\beta$ and $\theta_0=\pi$. The wavenumbers are $k=2$ (left) and $k=8$ (right). The resolution is intentionally kept low ($N_\rho=16$, $N_\theta=64$) to illustrate discrepancies. For small and moderate beam concentrations ($\beta=1,2,4$ at $k=2$, and $\beta=1,8$ at $k=8$) the numerical solution follows the analytic curve closely. At the most challenging case ($(k,\beta)=(8,64)$, green curve in the right figure), the resolution is insufficient, illustrating the onset of numerical difficulty at large~$\beta$ and higher frequency.}
\label{fig:gaussian-beam-scattering}
\end{figure}

In Fig.~\ref{fig:gaussian-beam-scattering}, we see the real part of the numerically computed scattered far field $u_\infty^+(\theta)$ with the explicit reference solution obtained from the modal spectrum \eqref{eq:VM-scattered-modes} for three values of the concentration parameter $\beta$ at two wavenumbers, $k=2$ and $k=8$, with $\theta_0=\pi$. 
Solid lines correspond to the analytic far field and open markers to the numerical solution evaluated on the outer boundary of the computational domain.  The resolution is deliberately kept low ($N_\theta=64$, $N_\rho=16$) so that discrepancies between the two solutions at the most challenging case are visible by eye.

For small and moderate beam concentrations the agreement is excellent.  Even at this relatively low resolution the numerical curves for $\beta=1,2,4$ when $k=2$ and for $\beta=1,8$ when $k=8$ are indistinguishable from their analytic counterparts.  The most demanding case in this experiment, $(k,\beta)=(8,64)$, is visibly under-resolved. This example illustrates the regime in which the method is challenged, namely when both the frequency $k$ and the concentration $\beta$ are large so that the effective angular bandwidth of the solution exceeds the available Fourier resolution and we get into the regime where the centrifugal barrier becomes relevant due to the excitation of high $m$-modes. Nonetheless, the solution can be accurately captured by moderately increasing the spectral resolution as demonstrated in the next section.

\subsubsection{Anisotropic media on the Poincaré--Beltrami disk}\label{sec:hyperbolic}

The numerical examples we discussed so far include radial profiles for the variable media and characteristic coordinates for the equation, which requires the eikonal solution. The advantages of this setup for illustration are the availability of the explicit solutions for radial media, and the simplicity of the equations in characteristic coordinates. However, the NIC-Helmholtz equation \eqref{eq:NIC-Helmholtz} presented in the main theorem is sufficiently general to handle non-radial cases and non-characteristic coordinates. In this section, we discuss the solution of more general scattering problems with anisotropic media using hyperbolic coordinates discussed in \eqref{eq:standard-hyperboloids}.



The anisotropic media example is derived from the refractive index \eqref{eqn:quadratic} by shifting its center. We place the center of the perturbation at the negative $x$-axis, \(\mathbf{x}_0=(-d,0)\). Using obstacle-centered polar coordinates \((r,\theta)\), we have $|\mathbf{x}-\mathbf{x}_0|^2 = r^2 + d^2 + 2 d r \cos\theta$. Thus, the off-centered refractive index can be written as
\begin{equation}\label{eq:anisotropic-refractive-index}
  n^2(r,\theta)
  = 1 + \frac{\kappa^2}{r^2 + d^2 + 2 d r \cos\theta + r_0^2},
\end{equation}
where \(r_0>0\) plays the role of a core radius: it regularizes the singularity at the center of the perturbation, $\{r=d, \theta=\pi\}$, and sets the maximum value of the refractive index. The medium is anisotropic since the refractive index depends on both \(r\) and \(\theta\).

Solving the eikonal equation for this anisotropic medium is possible, but nontrivial. A simpler approach is to use the standard hyperboloids in \eqref{eq:standard-hyperboloids} and the radial compactification \eqref{eq:alt-conformal}, thereby solving the NIC-Helmholtz equation on a punctured Poincaré-Beltrami disk. Specifically, we choose the spatial compactification
\begin{equation}\label{eq:hyperbolic-radius} 
  r = \frac{2\rho}{1-\rho^2}, \qquad G=\frac{1}{r'}=\frac{(1-\rho^2)^2}{2(1+\rho^2)}, 
\end{equation}
and the height function
\begin{equation}\label{eq:hyperbolic-height} 
  h=\pm \sqrt{1+r^2} = \pm \frac{1+\rho^2}{1-\rho^2}, \qquad H=\pm\frac{r}{\sqrt{1+r^2}} = \pm\frac{2\rho}{1+\rho^2}. 
\end{equation}
These choices satisfy \eqref{eq:compactification} and \eqref{eq:rate-condition}. Therefore, the rescaling by \eqref{eq:weight-rescaling-NIC} leads to the transformed equation \eqref{eq:NIC-Helmholtz} that is well-defined at the ideal boundary of the hyperbolic disk. In particular, the formally singular term in the compactified Helmholtz equation is regular on the entire $\rho$ domain including at \(\rho=1\). 
The NIC-Helmholtz equation reads:
\begin{align*}
 0 &= \partial_\rho \left[ \frac{(1-\rho^2)^2}{2(1+\rho^{2})}\,\partial_\rho u^\pm\right]
\pm \frac{4 i k \rho}{1+\rho^{2}}\,\partial_\rho u^\pm
+  \frac{1+\rho^{2}}{2\rho^{2}}\, \partial_\theta^2 u^\pm \\[6pt]
& + \Big[ k^2 \left(\frac{2}{1+\rho^2} + 
 \frac{1+\rho^2}{2} \frac{\kappa^2} {\rho^2+(d^2+r_0^2)\Omega^2+2 d \rho \Omega \cos\theta} \right) \pm \frac{2 i k (1-\rho^{2})}{(1+\rho^{2})^{2}}
+ \frac{1+\rho^{2}}{8\rho^{2}}
\Big] u^\pm.
\end{align*}
Another advantage of this compactification over the characteristic one is that the obstacle boundary condition is independent of the refractive index. In characteristic coordinates, the connection between the incident and scattered fields at the obstacle is established through the sound-soft condition that translates to \eqref{eq:boundary-variable}, which requires the eikonal solution. In hyperbolic coordinates, however, the time transformation is \emph{independent of the refractive index}. The relation between the incident and scattered fields at the obstacle is given by the height function, which satisfies asymptotic conditions \eqref{eq:hyperboloidal} but is otherwise free in the interior. In our case, the obstacle boundary condition for the rescaled scattered field is as follows:
\[
U^+|_\Gamma+U^-|_\Gamma=\frac{e^{ik|h|}}{\sqrt{r}} u^+|_\Gamma + \frac{e^{-ik|h|}}{\sqrt{r}} u^-|_\Gamma = 0 \quad \implies\quad u^+|_\Gamma = e^{-2ik |h|} u^-|_\Gamma.
\]
For scattering at the unit disk with $r=1$, the condition reads $u^+(|_\Gamma) = e^{-2\sqrt{2}ik} u^-|_\Gamma$. Note that the inner boundary with the compactification \eqref{eq:alt-conformal} is at \(\rho=\sqrt{2}-1\).



\begin{figure}[t]
  \centering
  \includegraphics[width=0.49\textwidth]{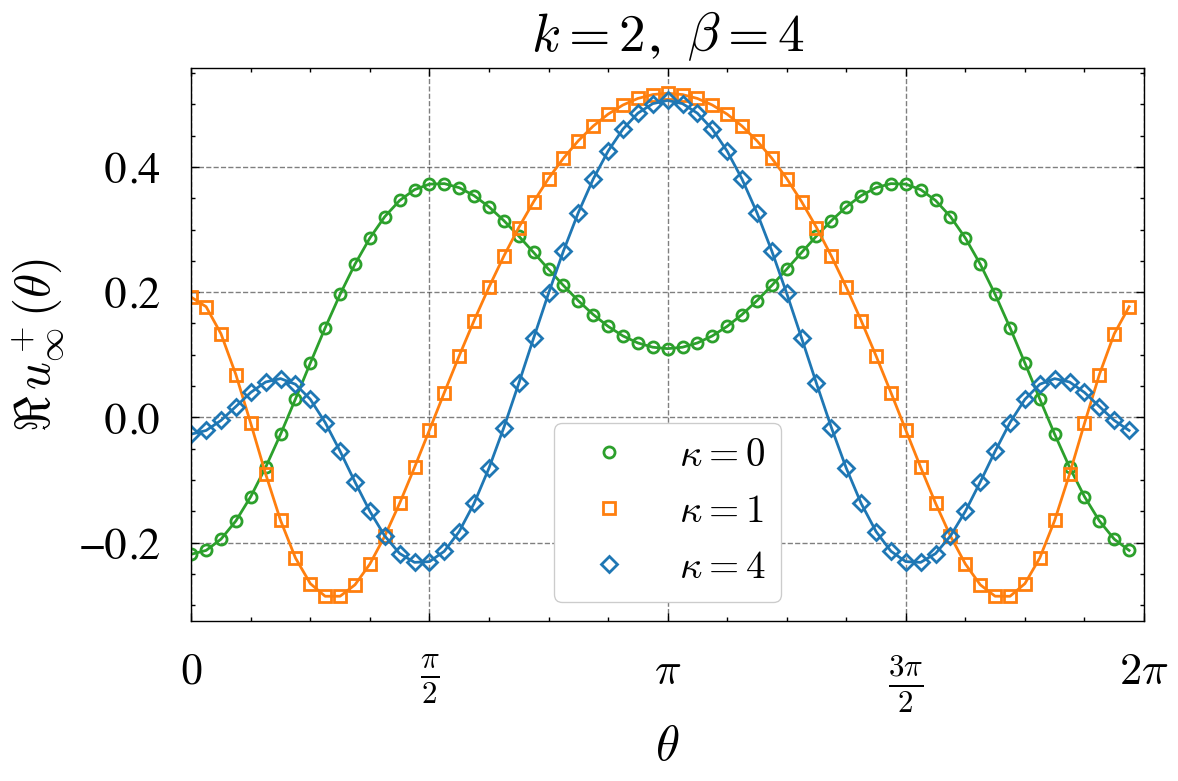}\hfill
  \includegraphics[width=0.49\textwidth]{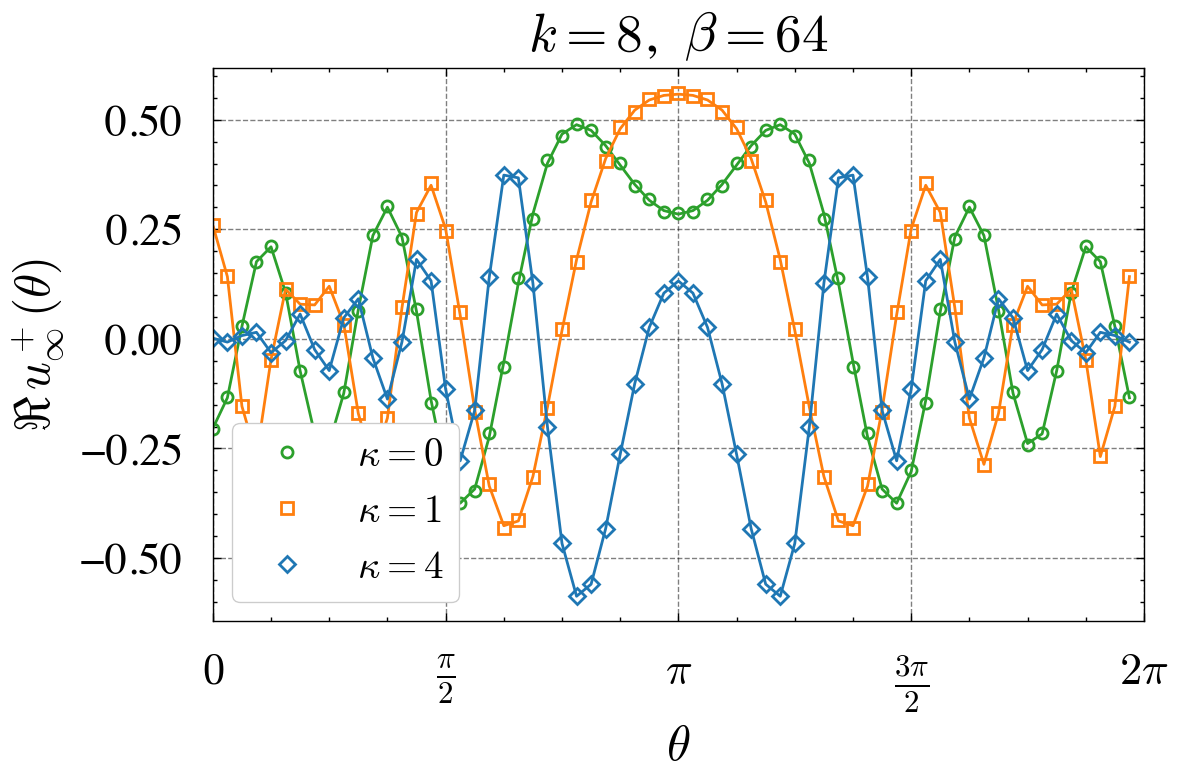}
\caption{Comparison of the exact scattered far field (solid lines) and the numerical solution (open markers) for Gaussian beam incidence \eqref{eq:gaussian} at different strenght parameters $\kappa$. We set $k=2$ (left) and $k=8$ (right) with $\beta=k^2$. The resolution is relatively high at $\{N_\rho=64$, $N_\theta=64\}$.}
\label{fig:gaussian-beam-scattering-quadratic}
\end{figure}

With $d=0$ and $r_0=0$, the anisotropic refractive index \eqref{eq:anisotropic-refractive-index} reduces to the radial index \eqref{eqn:quadratic} for which we have the explicit solution \eqref{eq:Sm-quadratic}. In Fig.~\ref{fig:gaussian-beam-scattering-quadratic}, we compare the numerically computed scattered far field $u_\infty^+(\theta)$ with the explicit reference solution. The reference solution for Gaussian beam incidence is constructed similarly to \eqref{eq:VM-scattered-modes} but with the scattering amplitudes replaced by the short-range solution from \eqref{eq:Sm-quadratic}. As we sum up all values of $m$ up to $N_\theta/2$, including ones with $|m|<\kappa k$, the modified mode number $\nu$ becomes complex. We plot the real part of the scattered far field for three values of the strength parameter, $\kappa=\{0, 1, 4\}$, at two wavenumbers $k=2$ (left) and $k=8$ (right) with concentration $\beta=k^2$ and $\theta_0=\pi$. The resolution is set to $\{N_\theta=64$, $N_\rho=64\}$ so the agreement between the numerical and explicit solutions is excellent for all cases. 

We include the $\kappa=0$ case in this experiment to show that the hyperbolic implementation matches the characteristic implementation at infinity (the green curves) even though the interior solution looks very different. The impact of increasing the strength of the medium through $\kappa$ is to change the relative phases of the angular modes, which modifies the far-field pattern through interference. The low-frequency case on the left panel of Fig.~\ref{fig:gaussian-beam-scattering-quadratic} shows a relatively simple pattern as only a few number of modes are excited. In the high-frequency strongly focused case on the right panel, many modes up to $|m|\sim k$ are excited and a larger fraction of them satisfy $|m|<\kappa k$, leading to a stronger intereference pattern at infinity. The NIC-Helmholtz implementation on the hyperbolic disk captures these effects accurately without solving the eikonal equation.


\begin{figure}[t]
  \centering
  \includegraphics[width=0.49\textwidth]{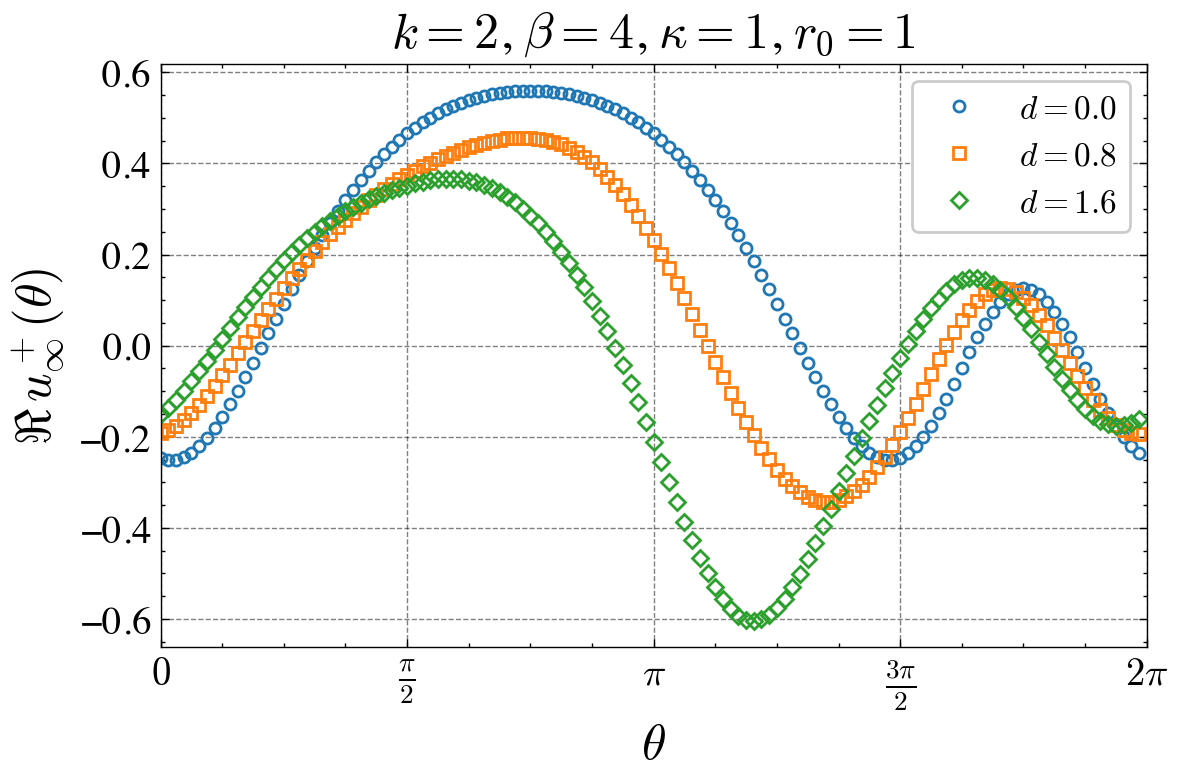}\hfill
  \includegraphics[width=0.49\textwidth]{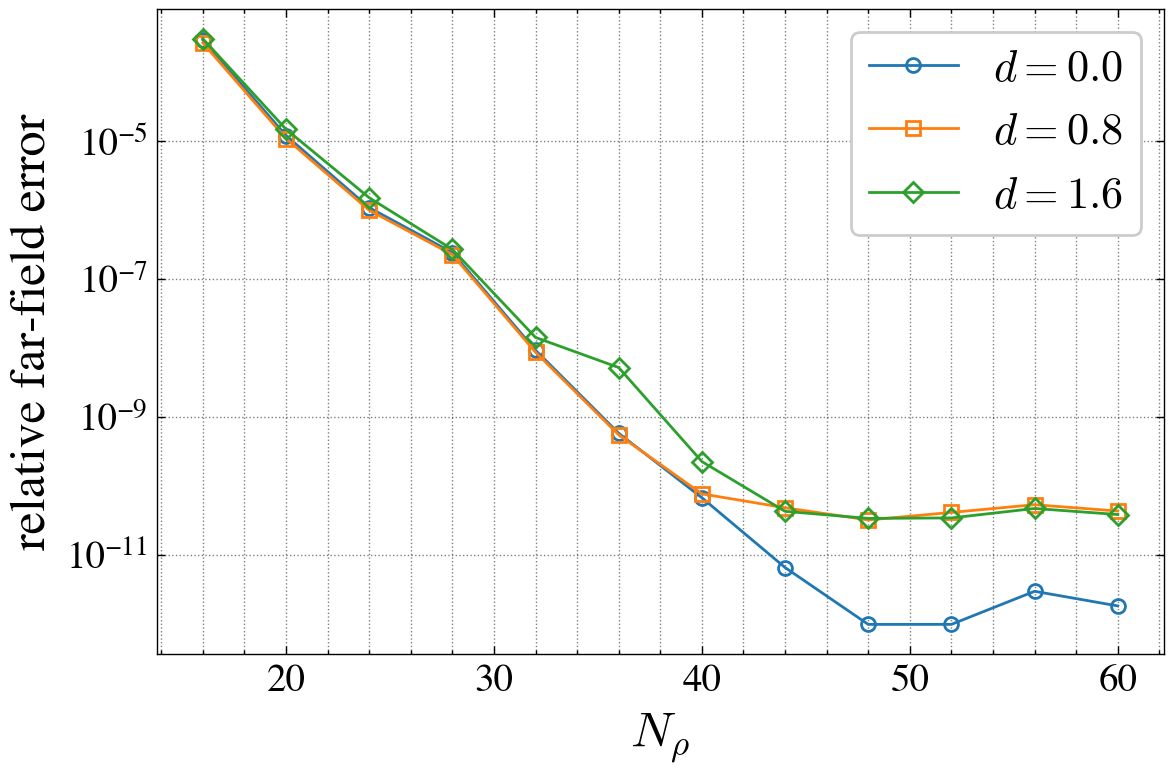}
\caption{Anisotropic medium scattering. Left: Real part of the scattered far field $u_\infty^+(\theta)$ for Gaussian beam incidence \eqref{eq:gaussian} with $k=2$, $\beta=4$, and $\theta_0=3\pi/4$. The refractive index is given by \eqref{eq:anisotropic-refractive-index} with $\kappa=r_0=1$ at different offset parameters $d$. Right: Self-convergence of the far field as the resolution $(N_\rho, N_\theta)$ is increased. The reference solution is computed at the highest resolution, $(N_\rho,N_\theta)=(64,128)$.}
\label{fig:anisotropic}
\end{figure}

Next, we solve the anisotropic case with the off--centered refractive index \eqref{eq:anisotropic-refractive-index}. Fixing $\kappa=1$ and $r_0=1$, the offset parameter $d$ controls the degree of angular variation in $n(\rho,\theta)$: the case $d=0$ corresponds to a radial medium; the case with $d>0$ leads to a fully non--separable, variable--coefficient, anisotropic operator in $(\rho,\theta)$. In Fig.~\ref{fig:anisotropic} we plot the real part of the numerically computed scattered far field $u_\infty^+(\theta)$
for a circular Gaussian beam with $k=2$, $\beta=4$, and $\theta_0=3\pi/4$. The incident beam is not aligned with the center of the medium and therefore the far field is not symmetric. As $d$ is increased from $d=0$ to $d=0.8$ to $d=1.6$, the side component becomes more pronounced. This behavior is due to the mode--mixing induced by anisotropy: unlike the radial case, the medium couples Fourier harmonics in $\theta$, and the far field arises from the interaction of many angular modes rather than from independent $m$--channels. 

Since no explicit scattering map is available in this non--radial setting, we confirm accuracy through a self--convergence study of the far field. On the right panel of Fig.~\ref{fig:anisotropic} we use the highest resolution as a reference and measure the relative $L^2$ error of $u_\infty^+(\theta)$ on $[0,2\pi)$. In all three cases, the error decays spectrally over several orders of magnitude as $N_\rho$ increases, reaching errors on the order of $10^{-11}$--$10^{-12}$. The isotropic baseline $d=0$ converges slightly further (down to $\sim 10^{-12}$) than the anisotropic cases, which level off around $\sim 10^{-11}$. These plots demonstrate that the NIC--Helmholtz formulation on the hyperbolic disk with a two--dimensional spectral discretization provides very accurate far-field computations even for challenging anisotropic media where separability is lost and neither the eikonal nor the explicit solutions are available.

\section{Discussion}\label{sec:discussion}

The core thesis of this paper is that \emph{infinity should be part of the computational domain} for the computation of scattering amplitudes and far fields. Even though there is no explicit notion of time in the frequency domain, Lorentzian geometry and time transformations provide a practical mechanism to incorporate infinity. The WKB-type rescaling \eqref{eq:scale_out} removes the asymptotic oscillatory decay from Helmholtz solutions, which leads to renormalized fields $u^\pm$ that are smooth up to the conformal boundary, $\mathscr I^\pm$.  In the resulting NIC formulation \eqref{eq:NIC-Helmholtz}, we can prescribe incoming radiation data $u^-_\infty$ at the past celestial sphere, $\mathscr I^-$, and extract outgoing data $u^+_\infty$ by local evaluation at the future celestial sphere, $\mathscr I^+$. Asymptotic observables such as far-field patterns and scattering amplitudes are obtained at the natural boundary where they are defined. This technique gives a numerical realization of the operator-theoretic picture of geometric scattering (Poisson operator and scattering matrix). 

Implementing this geometric framework numerically for variable media, we used a modified problem setup for the incident and scattered fields. In the standard decomposition, the incident field is typically a free-space plane wave, and all medium effects are captured by the scattered solution that solves an equation with a source term \eqref{eq:sourced-helmholtz} (the Lippmann--Schwinger viewpoint). The two-step solver introduced here defines incidence by mathematically well-posed incoming boundary data at $\mathscr I^-$ with finite energy, providing an incoming field that propagates in the background of the medium. The scattered field represents the obstacle correction. While this approach is computationally more intensive than the Lippmann--Schwinger approach, it has several advantages: (i) it cleanly separates the physics of incidence from the obstacle interaction; (ii) it avoids artificial truncation of the domain and associated absorbing boundary conditions or layers; and (iii) it is compatible with long-range media.

The numerical experiments presented in Sec.~\ref{sec:numerics} are primarily restricted to geometrically simple configurations (unit disk) and radial media in two space dimensions to enable direct validation against explicit solutions.  However, the method is quite general and can be applied to a wide range of scattering problems. In particular, we demonstrated in Sec.~\ref{sec:hyperbolic} that the NIC--Helmholtz equation can be solved accurately for anisotropic, unbounded media without the eikonal solution. 

There are several challenges and open questions for future work. The compactified NIC--Helmholtz operator is degenerate at the conformal boundary. It should be studied how this degeneracy affects the conditioning of discretizations. The current setup is restricted to asymptotically homogeneous media \eqref{eq:refractive}; more general anisotropic media will require compactifications and discretizations that preserve the regularity of $u^\pm$ up to $\mathscr I^\pm$. Numerical properties should be investigated in three dimensions, for general obstacles, and systems such as Maxwell. Such applications may connect this method with related work employing the semiclassical and microlocal perspective on scattering \cite{galkowski2022high, galkowski2024scattering, galkowski2025numerical} or computations using the Lippmann--Schwinger approach \cite{gillman2015spectrally, ying2015sparsifying}. 

Overall, NIC should be viewed as the beginning of a broader program: incorporating infinity as a geometric boundary aligns computation with scattering theory. I expect this perspective to be helpful in both analysis and algorithm design for wave propagation in unbounded domains.

\backmatter

\subsection*{Acknowledgements}

I thank Paul Martin, Mike O'Neil, and Eitan Tadmor for helpful discussions.
I used the software package Dedalus \cite{burns2020dedalus} for all 2D computations and ChatGPT 5.2 for coding (refactoring, plotting, translating). This material is supported by the National Science Foundation under Grant No. 2309084.

\medskip

\begin{appendices}

\section{Spacetime and geometric scattering}\label{sec:spacetime}
Scattering theory is, at heart, about how fields behave at infinity. The underlying geometric structure of spacetime is often hidden in signs, powers, and weights. To make this structure explicit, we first recall the global causal picture of Minkowski spacetime via Penrose's conformal compactification \cite{penrose_asymptotic_1963, penrose2011republication}. We then discuss how this picture is implicit in Melrose's scattering framework: the identification of incident and scattered fields in the compactified picture selects the boundaries $\scri^\pm$; the complex rescaling of variables corresponds to time transformations; and the natural half-density weight is connected to conformal properties of the Laplacian.

\subsection{The global causal structure of spacetime (from Penrose)}\label{sec:penrose}
A key idea in this section is that time transformations change the geometry of spatial slices. The level sets of standard time are Euclidean spaces. By transforming the time coordinate, we can get asymptotically hyperbolic or degenerate spaces. We exploit this change in asymptotic geometry to obtain conformally compact spatial slices of flat spacetime, leading to a regular radial compactification suitable for computation.

To demonstrate the global causal structure of spacetime relevant for numerical scattering, we consider the $(d+1)$-dimensional Minkowski spacetime $(\mathcal{M}, \eta)$, where $\mathcal{M}=\mathbb{R}_t\times \mathbb{R}_x^d$ and $\eta$ is the Lorentzian flat metric\footnote{For scattering, we are primarily interested in $d=2,3$}. In standard spherical coordinates $\{t,r,\omega\}$, with $r\ge 0$ and $\omega\in \mathbb{S}^{d-1}$, the metric takes the form
\[ \eta = -dt^2 + dr^2 + r^2 \gamma,\]
where $\gamma$ denotes the standard metric on the unit sphere $\mathbb{S}^{d-1}$. Radial compactification with $r=1/\rho$ results in a singular metric at infinity,
\begin{equation} \label{eq:singular}
\eta = -dt^2 + \frac{d\rho^2}{\rho^4} + \frac{\gamma}{\rho^2}.
\end{equation}
The singularity cannot be captured with a conformal rescaling, reflecting that Euclidean space is \emph{not} conformally compact.

\begin{figure}
  \centering
  \includegraphics[width=0.3\textwidth]{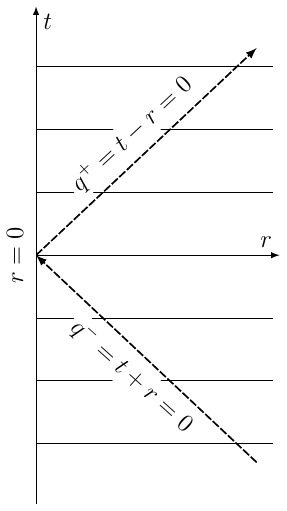}
    \hspace{1.5cm}
  \includegraphics[width=0.36\textwidth]{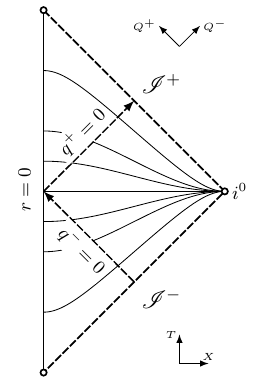}
  \caption{Left: Spacetime diagram depicting level sets of the standard time coordinate $t$ and the incoming/outgoing null cone at the origin, $q^\mp=0$ (each point of the diagram represents a sphere $\mathbb S^{d-1}$ of angles $\omega$ except for the line $r=0$). Right: Penrose diagram (see App.~\ref{app:penrose}) also depicting level sets of the standard time coordinate $t$ and null rays $q^\pm=0$. In the Penrose diagram, constant-$t$ slices meet at spatial infinity $i^0$ and therefore do not provide a regular description of the asymptotic behavior of radiative fields. Incoming and outgoing radiative fields live on $\mathscr I^-$ and $\mathscr I^+$, respectively.  
    }
    \label{fig:standard}
\end{figure}

Unlike the Euclidean case, the presence of a time direction endows the metric with a null cone at each point. We introduce null (or characteristic) coordinates,
\begin{equation} \label{eq:double_null}
q^+=t-r,\qquad q^-=t+r.
\end{equation}
The hypersurfaces $q^\pm=\text{const}$ are outgoing $(+)$ and incoming $(-)$ null cones. The set of null directions at a point forms a sphere $\mathbb{S}^{d-1}$. The Minkowski metric reads
\begin{equation}\label{eq:double-null-metric}
\eta=-(dq^\pm)^2 \mp 2\,dq^\pm\,dr+r^2\gamma.
\end{equation}
Radiative fields propagate along characteristics. To study their global behavior, Penrose introduced conformal compactification to relativity \cite{penrose_asymptotic_1963, penrose2011republication}, which maps infinity along null rays to finite endpoints at the conformal boundary, denoted $\scri^\pm$ (see App.~\ref{app:penrose} and Fig.~\ref{fig:standard}). In contrast to the Euclidean slicing in \eqref{eq:singular}, radial compactification along characteristics results in a conformally regular metric,
\begin{equation}\label{eq:conformal_null}
\rho^2 \eta = -\rho^2 (dq^\pm)^2\pm 2\,dq^\pm\,d\rho+\gamma.
\end{equation}
The generators $\partial_{q^\pm}$ are the natural time translations along $\scri^\pm$. Every cut of $\scri^\pm$ is a celestial sphere $\mathbb S^{d-1}$, and the null generators (integral curves of the null normal, which becomes tangent to $\scri^\pm$) foliate $\scri^\pm$ by these cuts. Radiative fields live on $\scri^\pm$.

The conformal viewpoint by Penrose matches the geometric scattering framework: incoming data are prescribed at past null infinity \(\scri^-\), while outgoing data are read off at future null infinity \(\scri^+\). Conformal compactification turns the limits \(r\to\infty\) and \(q^\pm\to\pm\infty\) into analysis \emph{at the boundary} \(\rho=0\), enabling stable local discretizations for unbounded-domain wave problems.

We generalize the transformation \eqref{eq:double_null} by using a height function, $h(r)$, with derivative $H(r):=\frac{dh}{dr}$ called boost, satisfying asymptotic conditions so that level sets of the new coordinate $\tau$ approach null coordinates asymptotically \cite{gowdy_wave_1981, friedrich1983cauchy, zenginouglu2008hyperboloidal},
\begin{equation} \label{eq:hyperboloidal}
\tau = t - h(r), \qquad H(r):=\frac{dh}{dr},\quad |H(r)| \leq 1 \ \text{for all }r, \ \text{and}\ \lim_{r\to\pm \infty} H(r) = \pm 1.
\end{equation}
The sign corresponds to the two null infinities, $\scri^\pm$. The transformed metric becomes
\[ \eta = -d\tau^2 - 2H d\tau dr + \left(1-H^2\right) dr^2 + r^2 \gamma. \]
Applying radial compactification $r=1/\rho$ and the conformal rescaling with $\rho^2$, we get
\begin{equation}\label{eq:conformal_hyperboloidal}
\rho^2\eta = -\rho^2d\tau^2 + 2H d\tau d\rho + \frac{1-H^2}{\rho^2} d\rho^2 + \gamma.
\end{equation}
To obtain a smooth conformal limit at infinity, it is convenient to assume the following asymptotics for the boost
\begin{equation}\label{eq:hyperboloidal-nondeg}
1-H(r)^2 = \frac{a^2}{r^2} + o(r^{-2})
\quad\Longleftrightarrow\quad
H(r) = \pm\left(1 - \frac{a^2}{2r^2} + o(r^{-2})\right),
\end{equation}
so that $\rho^2 \eta_\tau \to a^2\,d\rho^2 + \gamma$ as $\rho\to 0$. In asymptotically null coordinates (with $a=0$) \cite{zenginouglu2024hyperbolic}, the limiting metric degenerates. In null coordinates ($H=\pm 1$), we obtain \eqref{eq:conformal_null}. The conditions in the main theorem, \eqref{eq:compactification} and \eqref{eq:rate-condition} ensure that the transformed spacetime metric is conformally compact with a smooth extension to the conformal boundary. Below is a non-characteristic example that satisfies the above conditions.

\begin{example}[Standard hyperboloids]\label{eq:example}
The level sets of $\tau = t-\sqrt{1+r^2}$ are hyperboloids with
\begin{equation}\label{eq:standard-hyperboloids}
h(r)=\sqrt{1+r^2}
= r + \frac{1}{2r} + \mathcal{O}(r^{-3}),
\qquad
H(r)=\frac{r}{\sqrt{1+r^2}}
= 1 - \frac{1}{2r^2}+\mathcal{O}(r^{-4}),
\end{equation}
The height function satisfies \eqref{eq:hyperboloidal} and \eqref{eq:hyperboloidal-nondeg}. The induced geometry is the hyperbolic ball with constant negative sectional curvature. To see this, we use a slightly different compactification
\begin{equation} \label{eq:alt-conformal} r = \frac{2\rho}{1-\rho^2}. \end{equation}
The induced metric is then precisely the metric on the Poincaré--Beltrami model of hyperbolic space,
\[
\eta_\tau = g_{\mathbb{H}^2} = \frac{4}{(1-\rho^2)^2}\left(d\rho^2 + \rho^2 \gamma\right).
\]
\end{example}

The discussion above demonstrates that we can control the geometry of the spatial slices by choosing appropriate asymptotic conditions on the time transformation \eqref{eq:hyperboloidal}. The spatial metric on the standard $t=\mathrm{const}$ surfaces is Euclidan. By a suitable choice of height function, we can obtain asymptotically hyperbolic or degenerate spaces, both of which are conformally compact and allow for radial compactification suitable for computation. The sign of the height function determines whether the slices approach past or future null infinity.

In summary, the construction of a regular compactification at null infinity consists of three steps: (i) radial compactification, (ii) time transformation, (iii) and a conformal rescaling. In the next subsection, we discuss how each of these steps can be incorporated into geometric scattering theory. The resulting partial differential equations are regular at the boundary, thereby allowing us to compute numerically the scattering amplitudes on the celestial sphere at null infinity.

\subsection{Geometric scattering theory (to Melrose)}\label{sec:melrose}
To make infinity part of the computational formalism, we connect the conformal treatment of null infinity in spacetime by Penrose \cite{penrose_asymptotic_1963, penrose2011republication}) and the operator-theoretic framework of scattering resonances on noncompact spaces by Melrose \cite{melrose1995geometric, melrose2020spectral, dyatlov2019mathematical}. Here, I summarize the key points heuristically without technical proofs. The bridge connecting these two topics is organized by the three steps discussed above: radial compactification in Sec.~\ref{sec:radial}, time transformation in Sec.~\ref{sec:time}, and conformal rescaling in Sec.~\ref{sec:conformal}.

\subsubsection{Radial compactification and Melrose's scattering calculus}\label{sec:radial}
In his treatment of conformal infinity, Penrose performs a compactification along null directions, thereby replacing asymptotic limits by local analysis at the conformal boundary. Similarly, Melrose treats infinity as a microlocal boundary on which one can do analysis. For $\mathbb{R}^d$, radial compactification produces a compact manifold $X$ with boundary $\partial X \simeq \mathbb{S}^{d-1}$ and an asymptotic boundary-defining function $\rho = r^{-1}$, so hat $\rho\to 0$ corresponds to $r\to\infty$ \cite{melrose1995geometric, melrose2020spectral, vasy1998geometric}. The Euclidean metric takes the form of a \emph{scattering metric} (compare \eqref{eq:singular}),
\[
  g \sim \frac{d\rho^2}{\rho^4} + \frac{\gamma}{\rho^2}
  \quad\text{near }\partial X,
\]
Unlike the conformally compact case \eqref{eq:conformal_null}, the normal component is strongly degenerate, $\rho^{-4} d\rho^2$, so $g$ does not extend to $\partial X$ as a smooth metric on $TX$ but on the b-tangent bundle ${}^{\mathrm b}TX$. The associated Helmholtz operator is \emph{singular} at $\rho=0$ when viewed on the ordinary tangent/cotangent bundles (see \eqref{eq:singular-helmholtz}). 
Melrose's resolution is to replace the ordinary geometry by scattering geometry at $\partial X$. The Lie algebra of vector fields of uniformly bounded length with respect to a scattering metric is
\[
  \mathcal{V}_{\mathrm{sc}}(X)=x\,\mathcal{V}_{\mathrm b}(X),
  \quad\text{locally spanned by}\quad \{x^2\partial_x,\ x\partial_{y_1},\dots,x\partial_{y_{d-1}}\},
\]
which is the space of smooth sections of the scattering tangent bundle ${}^{\mathrm{sc}}TX$ and dual bundle ${}^{\mathrm{sc}}T^*X$ \cite{melrose1995geometric}. In this geometry, $\Delta_g\in\mathrm{Diff}^2_{\mathrm{sc}}(X)$ is a uniformly degenerate second-order operator which is elliptic with respect to the scattering principal symbol. The scattering Sobolev scales $H^{s,\alpha}_{\mathrm{sc}}$ encode decay (weight) and regularity uniformly up to $\partial X$.
Radiation conditions are then encoded microlocally at infinity (at the radial sets) via the choice of weights and propagation direction (outgoing vs incoming resolvents), rather than as a naive local boundary condition at $\rho=0$.
This framework also admits long-range perturbations (metrics and potentials) with controlled polyhomogeneous behavior at $\partial X$
\cite{vasy1998geometric}.

While Melrose's scattering framework rigorously localizes the problem at infinity, its direct translation into a numerical discretization is not obvious. By contrast, time transformations adopted to characteristics lead to conformally compact spatial slices (Sec.~\ref{sec:penrose}). The challenge is to incorporate this spacetime regularization into time-harmonic scattering, which lacks a notion of time.



\subsubsection{Time transformation}\label{sec:time}
The connection between time transformation and frequency-domain variables is obtained through the time-harmonic assumption. A time-harmonic field has the form $V(t,x)=e^{-ikt}\,U(x)$. Applying the time transformation \eqref{eq:hyperboloidal}, we get
\[
V(\tau+h(x),x)= e^{-ik(\tau+h(x))}\,U(x)
= e^{-ik\tau}\,\left(e^{-ikh(x)}U(x)\right).
\]
Thus, the time-shifted scalar field is time-harmonic at frequency $k$. By defining the rescaled field in the parenthesis as $u(x)$ and applying the inverse transform, we observe that the time shift corresponds to multiplication by $e^{ikh}$ via $U(x)=e^{ikh} u(x)$. The rescaling induces a space-dependent phase shift on the amplitude \cite{zenginouglu2011geometric}. To see the effect of the time transformation on the Helmholtz operator, let \(L(k):=\Delta+k^2 \).  Conjugating \(L(k)\) by the phase \(e^{ik h}\), using the identity \(e^{-ik h}\,\nabla\,e^{ik h}=\nabla+i k\,\nabla h\), we get
\begin{equation}\label{eq:conj-helm}
L_h(k) = e^{-ik h}\,L(k)\,e^{ik h} =\Delta + 2ik\,\nabla h\cdot\nabla + ik\,\Delta h + k^2\left(1-|\nabla h|^2\right).
\end{equation}
Assuming a radial height function $h(x)=h(r)$, the asymptotic condition \eqref{eq:hyperboloidal}  implies $|\nabla h|\to 1$ as $r\to\infty$. Then, due to our assumption \eqref{eq:refractive} on the refractive index, the \(k^2\)-term approaches zero at \(\mathscr I^\pm\). The time transformation asymptotically flattens the Helmholtz amplitude by removing the leading  WKB oscillations from $U$ and yields a more slowly varying amplitude $u$ \cite{zenginouglu2021null}. 

Equivalently, the map $U\mapsto u=e^{-ik h}U$ is a semiclassical gauge transform (with semiclassical parameter $1/k$) that shifts covectors by $-\nabla h$.
With $h$ chosen so that $|\nabla h|\approx 1$ asymptotically, the transformed operator \eqref{eq:conj-helm} is aligned with outgoing characteristics,
and after compactification the boundary at infinity becomes characteristic: the leading asymptotic behavior is governed by the transport term
$2ik\,\nabla h\cdot\nabla$ rather than by a nonzero $k^2$ potential.

However, the time transformation is not sufficient to obtain a regular compactification at infinity. Continuing with a radial height function $h=h(r)$ and writing $H:=h'(r)$, the conjugated free-space operator $L_h(k)=e^{-ikh}(\Delta+k^2)e^{ikh}$ becomes
\begin{equation}\label{eq:Lhk}
L_{h}(k) =\partial_r^2+\frac{d-1}{r}\partial_r+\frac{1}{r^2}\Delta_{\mathbb S^{d-1}} +2ik\,H \partial_r +ik\left(\partial_r H+\frac{d-1}{r}H\right) +k^2\left(1-H^2\right).
\end{equation}
Under radial compactification $\rho=1/r$ (and extracting the natural factor $\rho^2$) we obtain
\begin{align*}
\rho^{-2} L_{h}(k)
=&(\rho\,\partial_\rho)^2+(2-d)(\rho\,\partial_\rho)+\Delta_{\mathbb S^{d-1}}
-2ik\,H \partial_\rho\\
&+ik \left(-\partial_\rho H+\frac{d-1}{\rho}H \right)
+k^2\frac{1-H^2}{\rho^2},
\end{align*}
With \eqref{eq:hyperboloidal-nondeg} implying \(H=1-a^2\rho^2\), the term
$k^2(1-H^2)/\rho^2$ extends smoothly to $\rho=0$. However, the lower-order term $\frac{d-1}{\rho}H$ diverges like $(d-1)/\rho$ for $d\ge2$. Thus, demodulation by $e^{ikh}$ removes the oscillatory phase but does not by itself produce a smooth compactified operator; an additional rescaling is required.


\subsubsection{Conformal rescaling}\label{sec:conformal}
The final ingredient for a regular compactification of the Helmholtz operator at infinity is the conformal rescaling. This rescaling can be understood in many equivalent ways: physically, conformally, algebraically, and microlocally.

Physically, the rescaling is related to flux conservation and the geometric divergence of spherical surfaces towards infinity. The surface area of spheres in $d$-dimensional space grows as $r^{d-1}$. As the energy of an outgoing wave is distributed on the surface of the sphere, its amplitude decays like $r^{(d-1)/2}$. The rescaling takes out this decay.

Conformally, the rescaling arises from the transformation properties of derivative operators. We recap the standard calculation \cite{penrose2011republication, wald1984general} using the conjugation of the conformally covariant Yamabe operator \cite{yamabe1960deformation}. For an $n$-dimensional manifold $(\mathcal{M},g)$, the Yamabe operator is defined as
\[
\mathcal{Y}_g := \nabla^2_g - \frac{n-2}{4(n-1)} R[g]
\]
With respect to a conformal metric given via $g = \Omega^{-2}\tilde g$, the Yamabe operator acting on scalars satisfies
\begin{equation}
\label{eq:yamabe}
\Omega^{\frac{n-2}{2}} \,\mathcal{Y}_g \,\Omega^{-\frac{n-2}{2}}=\Omega^{-2} \mathcal{Y}_{\tilde{g}}
\end{equation}
Now consider the manifold to be Lorentzian with spacetime dimension $n=d+1$. Given a solution to the scalar wave equation $\nabla^2_g V=0$, we obtain from the transformation of the Yamabe operator,
\begin{equation}
\label{eq:conf-coupled}
\left(\,\nabla^2_{\tilde g} - \frac{d-1}{4d}\,\tilde R\,\right)\tilde V = 0, \qquad \tilde V=\Omega^{-\frac{d-1}{2}}V.
\end{equation}
The coefficients of this conformal equation are smooth, including at null infinity, for a conformally regular metric. In Eqs.~\eqref{eq:conformal_null} and \eqref{eq:conformal_hyperboloidal}, we used the conformal factor $\Omega=\rho=1/r$ and obtained a regular metric at null infinity. The conformal factor may take a different coordinate expression (for example, \eqref{eq:alt-conformal}) but decays as $\Omega\sim 1/r$. The rescaling of the variable $V$ in \eqref{eq:conf-coupled} counteracts the decay of the field such that the rescaled field has a finite, typically nonzero limit. Thus, $\tilde V$ extends to $\scri$ and satisfies the regular equation \eqref{eq:conf-coupled} on the conformal completion. 

Algebraically, the rescaling is known as the Liouville transform and is widely used in the analysis of wave equations. The Laplacian satisfies,
\[
r^{\frac{d-1}{2}}\left(\partial_r^2+\frac{d-1}{r}\partial_r+\frac{1}{r^2}\Delta_{\mathbb S^{d-1}}\right)r^{-\frac{d-1}{2}}
=\partial_r^2+\frac{1}{r^2}\Delta_{\mathbb S^{d-1}}
+\frac{(d-1)(3-d)}{4r^2}.
\]
Then, the conjugated operator $\widetilde L_{h}(k):= r^{\frac{d-1}{2}}\, L_{h}(k)\,r^{-\frac{d-1}{2}}$ reads
\begin{equation}\label{eq:Lhk_rad}
\widetilde L_h(k)=\partial_r^2+\frac{1}{r^2}\Delta_{\mathbb S^{d-1}}
+\frac{(d-1)(3-d)}{4r^2}
+2ik\,H \partial_r
+ik\, \partial_r H
+k^2\left(1-H^2\right).
\end{equation}
Comparing with \eqref{eq:Lhk}, we see that the $1/r$ terms drop out. In particular, there is no longer a singular \(\frac{H}{r}\) term. Upon radial compactification, we get
\begin{equation}\label{eq:tLhk}
\rho^{-2} \widetilde L_h(k)= \rho^2\partial_\rho^2
+ 2\left(\rho-ik \,H\right)\partial_\rho
+ \Delta_{\mathbb S^{d-1}}
+ \frac{(d-1)(3-d)}{4}
- ik\,\partial_\rho H
+ k^2\,\frac{1 - H^2}{\rho^2}\,.
\end{equation}
If the boost $H$ satisfies the asymptotic condition \eqref{eq:hyperboloidal-nondeg}, there are no divergent terms in the equation. 

Microlocally, the rescaling implements the half-density normalization, putting far-field amplitudes in \(L^2(\mathbb S^{d-1})\). 
Melrose formulated the scattering calculus on half-densities to obtain coordinate-invariant principal symbols and adjoint operations~\cite{melrose1995geometric}. 
Under a change of local coordinates \(x \mapsto \tilde{x}\), a half density transforms by the square root of the Jacobian determinant, $|\det(D\tilde{x}/Dx)|^{1/2}$,
so that the tensor product of two half-densities defines an invariant density. 
In spherical coordinates on \(\mathbb{R}^d\), the Euclidean density is \( |dx| = r^{d-1} \, dr \, |d\omega| \), 
hence a half-density carries the factor \( |dx|^{1/2} = r^{\frac{d-1}{2}} (dr \, |d\omega|)^{1/2} \). 
This factor gives rise to the familiar \(r^{(d-1)/2}\) Liouville or conformal scaling, which naturally places far-field amplitudes in \(L^2(\mathbb S^{d-1})\).

These different viewpoints, physical, conformal, algebraic, and microlocal, are different faces of the Yamabe identity \eqref{eq:yamabe}. The conformal picture is particularly powerful because it allows us to choose a more general scaling or compactification compatible with the asymptotic behavior, while still obtaining regular equations at null infinity. This freedom is used, for example, in the construction of the null infinity layer in \cite{zenginouglu2021null}, and the computation on the hyperbolic plane in Sec.~\ref{sec:hyperbolic}.


\section{Penrose diagram}\label{app:penrose}

An illustrative way to represent the global causal structure of spacetimes is via a \emph{Penrose diagram}, obtained by conformal compactification \cite{penrose_asymptotic_1963, penrose2011republication}. Following Penrose, we compactify along the null directions by applying the tangent map to the null coordinates:
\[
Q^+=\tan^{-1}q^+,\qquad Q^-=\tan^{-1}q^-,\qquad Q^\pm\in\left(-\frac{\pi}{2},\frac{\pi}{2}\right).
\]
The Minkowski metric in $(q^\pm,\theta)$ and $(Q^\pm,\theta)$ coordinates reads
\begin{align*}\label{eq:double_null_metric}
\eta&=-\,dq^+\,dq^-+\frac{(q^--q^+)^2}{4}\,d\theta^2\\
&=\frac{1}{\cos^2Q^+\,\cos^2Q^-}\left(-\,dQ^+\,dQ^-+\frac{1}{4}\sin^2(Q^--Q^+)\,d\theta^2\right).
\end{align*}
The overall factor $\Omega=\cos Q^+\,\cos Q^-$ is the \emph{conformal factor} capturing the blow-up as $Q^\pm\to\pm\pi/2$. Working with the rescaled metric $\bar{\eta}=\Omega^2 \eta$ makes the geometry regular at the boundary. Introducing Penrose time $T=Q^++Q^-$ and radius $X=Q^--Q^+$ with $X\in[0,\pi)$ (the restriction $r\ge 0$ implies $Q^-\ge Q^+$, hence $X\ge0$), 
the physical image of Minkowski space becomes the open region
$-\pi+X<T<\pi-X$ with $X\in[0,\pi)$. The conformal boundary is the hypersurface $\Omega=0$, or equivalently $|T|+X=\pi$. The conformal factor $\Omega$ acts as the boundary-defining function for the conformal completion of flat spacetime.

Of particular relevance for frequency-domain scattering is the conformal boundary along null geodesics. \emph{Future null infinity} is
$\scri^+=\{Q^-=\frac{\pi}{2},-\frac{\pi}{2}<Q^+<\frac{\pi}{2}\},$
and \emph{past null infinity} is $\scri^-=\{Q^+=-\frac{\pi}{2},-\frac{\pi}{2}<Q^-<\frac{\pi}{2}\}$. Each is a connected null cylinder with topology $\mathbb{R}\times S^1$: the $\mathbb{R}$ factor is the relevant null time ($Q^+$ on $\scri^+$ or $Q^-$ on $\scri^-$), while the $S^1$ is the space of asymptotic angles $\theta$. In a 2D Penrose diagram where the $S^1$ is suppressed, $\scri^\pm$ appears as the two diagonal edges (see Fig.~\ref{fig:standard}).



\end{appendices}

\bibliography{refs}

\end{document}